\documentclass{amsart}
\usepackage{amssymb,latexsym}
\usepackage{amsfonts}

\newtheorem{thm}{Theorem}[section]
\newtheorem{lem}[thm]{Lemma}

\newtheorem{rem}[thm]{Remark}
\newtheorem{pro}[thm]{Proposition}
\newtheorem{defi}[thm]{Definition}

\newcommand{\ba}{\begin{array}}
\newcommand{\ea}{\end{array}}

\def \qed{\cqfd}

\newcommand*{\QEDB}{\hfill\ensuremath{\square}}
\def\qed{\vbox{\hrule
\hbox{\vrule\hbox to 5pt{\vbox to 8pt{\vfil}\hfil}\vrule}\hrule}}

\newcommand{\beg}{\begin{eqnarray*}}
\newcommand{\begn}{\begin{eqnarray}}
\newcommand{\en}{\end{eqnarray*}}
\newcommand{\enn}{\end{eqnarray}}

\begin{document}
\title[The conical K\"ahler-Ricci flow with weak initial data]{The conical K\"ahler-Ricci flow with weak initial data on Fano manifold}
\author{JiaWei Liu}
\address{Jiawei Liu\\Beijing International Center for Mathematical Research\\ Peking University \\ Beijing 100871\\ China\\} \email{jwliu@math.pku.edu.cn}
\author{Xi Zhang}
\address{Xi Zhang\\Key Laboratory of Wu Wen-Tsun Mathematics\\ Chinese Academy of Sciences\\School of Mathematical Sciences\\
University of Science and Technology of China\\
Hefei, 230026, China\\ } \email{mathzx@ustc.edu.cn}
\thanks{}

\begin{abstract} In this paper, we prove the long-time existence and uniqueness of the conical K\"ahler-Ricci flow with weak initial data which admits $L^{p}$ density for some $p>1$ on Fano manifold. Furthermore, we  study the convergence behavior of this flow.

\end{abstract}

\maketitle
\section{Introduction}
\setcounter{equation}{0}

Conical K\"ahler-Einstein metric plays an important role in solving the Yau-Tian-Donaldson's conjecture (see \cite{CDS1, CDS2, CDS3, T1}). There has been renewed interest in conical K\"ahler-Einstein metric in recent years, see references \cite{RB, SB, CGP, EGZ, GP1, JMR, LS, SW}etc. On the other hand, the conical K\"ahler-Ricci flow was introduced to attack the existence problem of  conical K\"ahler-Einstein metric.
 The long-time existence and limit behaviour of the conical K\"ahler-Ricci flow has been widely studied. Chen-Wang \cite{CW} studied the strong conical K\"ahler-Ricci flow and obtained the short-time existence, Y.Q. Wang \cite{YQW} and the authors \cite{JWLXZ} got the long-time existence of the conical K\"ahler-Ricci flow respectively. In \cite{JWLXZ}, the authors also considered the convergence of this flow on Fano manifold with positive twisted first Chern class, they proved that, for any cone angle  $0<2\pi \beta < 2\pi $, the conical K\"ahler-Ricci flow converges to a conical K\"ahler-Einstein metric if there exists one. Chen-Wang \cite{CW1} obtained the convergence result of this flow when the twisted first Chern class is negative or zero. Later, by adopting the smooth approximation used by the authors in \cite{JWLXZ}, L.M. Shen \cite{LMSH1}\cite{LMSH2} generalized Song-Tian's result \cite{JSGT} and Tian-Zhang's result \cite{GTZZ} to the unnormalized conical K\"ahler-Ricci flow,  and G. Edwards \cite{GEDWA} obtained the uniform scalar curvature bound when the twisted first Chern class is negative on the foundation of Song-Tian's result \cite{JSGT1}.

In \cite{JWLXZ}, the authors studied the conical K\"ahler-Ricci flow which starts with a model conical K\"ahler metric
\begin{eqnarray}\label{01}
\omega_{\beta}=\omega_{0}+\sqrt{-1}k\partial\bar{\partial}|s|_{h}^{2\beta}
\end{eqnarray}
on Fano manifold, where $\omega_{0}\in c_{1}(M)$ is a smooth K\"ahler metric, $s$ is the defining section of a smooth divisor $D\in|-\lambda K_{M}|$ and $h$ is a smooth Hermitian metric on $-\lambda K_{M}$ with curvature $\lambda\omega_{0}$.  In \cite{CW, CW1, YQW}, Chen-Wang studied the existence of the conical K\"ahler-Ricci flow from initial $(\alpha,\beta)$ metric or weak $(\alpha,\beta)$ metric with other assumptions.

In this paper, we mainly study the long-time existence, uniqueness and convergence of the conical K\"ahler-Ricci flow with some weak initial data which admits $L^{p}$ density with $p>1$ on Fano manifold. We consider the conical K\"ahler-Ricci flow by using smooth approximation of the twisted K\"ahler-Ricci flows as that in \cite{JWLXZ}.

  Let $M$ be a Fano manifold with complex dimension $n$, $\omega_{0}\in c_{1}(M)$ be a smooth K\"ahler metric. For any $p\in(0,\infty]$, we define the class
\begin{eqnarray}\label{02}
\mathcal{E}_{p}(M,\omega_{0})=\big\{\varphi\in\mathcal{E}(M,\omega_{0}) | \frac{(\omega_{0}+\sqrt{-1}\partial\bar{\partial}\varphi)^{n}}{\omega_{0}^{n}}\in L^{p}(M,\omega_{0}^{n}) \big\},
\end{eqnarray}
where the class
\begin{eqnarray}\label{03}
\mathcal{E}(M,\omega_{0})=\big\{\varphi\in PSH(M,\omega_{0}) | \int_{M} (\omega_{0}+\sqrt{-1}\partial\bar{\partial}\varphi)^{n}= \int_{M} \omega^{n}_{0}\big\}
\end{eqnarray}
defined by Guedj-Zeriahi in \cite{VGAZ} is the largest subclass of $PSH(M,\omega_{0})$ on which the operator $(\omega_{0}+\sqrt{-1}\partial\bar{\partial}\cdot)^{n}$ is well defined and the comparison principle is valid.
When $p>1$, by  S. Kolodziej's $L^{p}$-estimate \cite{K} and S. Dinew's uniqueness theorem \cite{SDIN} (see also Theorem $B$ in \cite{VGAZ}), we know that the functions in $\mathcal{E}_{p}(M,\omega_{0})$ are H\"older continuous with respect to $\omega_{0}$ on $M$.

Let $D$ be a divisor on $M$. By saying a closed positive $(1,1)$-current $\omega\in2\pi c_{1}(M)$ with locally bounded potential is a conical K\"ahler metric with cone angle $2\pi \beta $ ($0<\beta < 1$ ) along $D$, we mean that $\omega$ is a smooth K\"ahler metric on $M\setminus D$, and near each point $p\in D$, there exists local holomorphic coordinate $(z^{1}, \cdots, z^{n})$ in a neighborhood $U$ of $p$ such that locally $D=\{z^{n}=0\}$, and $\omega$ is asymptotically equivalent to the model conical metric
\begin{eqnarray}
\sqrt{-1} |z^{n}|^{2\beta -2} dz^{n}\wedge d\overline{z}^{n} +\sqrt{-1} \sum_{j=1}^{n-1} dz^{j}\wedge d\overline{z}^{j}\ \ \ \ \ \ \ \ on\ \ \ U.
\end{eqnarray}

Assume that $D\in|-\lambda K_{M}|$ $(\lambda\in\mathbb{Q})$, $\mu=1-(1-\beta)\lambda$, $\hat{\omega}\in c_{1}(M)$ is a K\"ahler current which admits $L^{p}$ density with respect to $\omega_{0}^{n}$ for some $p>1$ and $\int_{M}\hat{\omega}^{n}=\int_{M}\omega_{0}^{n}$. We study the long-time existence, uniqueness and convergence of the following conical K\"ahler-Ricci flow with weak initial data $\hat{\omega}$
\begin{eqnarray}\label{TCKRF1}
\begin{cases}
  \frac{\partial \omega(t)}{\partial t}=-Ric(\omega(t))+\mu\omega(t)+(1-\beta)[D].\\
  \\
  \omega(t)|_{t=0}=\hat{\omega}\\
  \end{cases}
  \end{eqnarray}

From now on, we denote the K\"ahler current $\hat{\omega}=\omega_{0}+\sqrt{-1}\partial\bar{\partial}\varphi_{0}:=\omega_{\varphi_{0}}$ with $\varphi_{0}\in\mathcal{E}_{p}(M,\omega_{0})$ for some $p>1$.
\begin{defi}\label{04.5}We call $\omega(t)$ a long-time solution to the conical K\"ahler-Ricci flow $(\ref{TCKRF1})$ if it satisfies the following conditions.
\begin{itemize}
  \item  For any $[\delta, T]$ ($\delta , T>0$), there exist constant $C$ such that
   \begin{eqnarray*}
  C^{-1}\omega_{\beta}\leq\omega(t)\leq C\omega_\beta\ \ \ \ \ \ \ \ \ on\ \ [\delta, T]\times (M\setminus D);
  \end{eqnarray*}
  \item  On $(0,\infty)\times(M\setminus D)$, $\omega(t)$ satisfies the smooth K\"ahler-Ricci flow;
  \item On $(0,\infty)\times M$, $\omega(t)$ satisfies equation (\ref{TCKRF1}) in the sense of currents;
  \item There exists metric potential $\varphi(t)\in C^{0}\big([0,\infty)\times M\big)\cap C^{\infty}\big((0,\infty)\times (M\setminus D)\big)$ such that $\omega(t)=\omega_{0}+\sqrt{-1}\partial\bar{\partial}\varphi(t)$ and $\lim\limits_{t\rightarrow0^{+}}\|\varphi(t)-\varphi_{0}\|_{L^{\infty}(M)}=0$;
  \item On $[\delta, T]$, there exist constant $\alpha\in(0,1)$ and $C^{\ast}$ such that the above metric petential $\varphi(t)$ is $C^{\alpha}$ on $M$ with respect to $\omega_{0}$ and $\| \frac{\partial\varphi(t)}{\partial t}\|_{L^{\infty}(M\setminus D)}\leqslant C^{\ast}$.
  \end{itemize}
\end{defi}
In Definition \ref{04.5}, by saying $\omega(t)$ satisfies equation (\ref{TCKRF1}) in the sense of currents on $(0,\infty)\times M:=M_{\infty}$, we mean that for any smooth $(n-1,n-1)$-form $\eta(t)$ with compact support in $(0,\infty)\times M$, we have
  \begin{eqnarray*}
  \int_{M_{\infty}}\frac{\partial \omega(t)}{\partial t}\wedge \eta(t,x)dt=\int_{M_{\infty}}(-Ric(\omega(t))+\mu\omega(t)+(1-\beta)[D])\wedge \eta(t,x)dt,
  \end{eqnarray*}
  where the integral on the left side can be written as
    \begin{eqnarray*}
  \int_{M_{\infty}}\frac{\partial \omega(t)}{\partial t}\wedge \eta(t,x)dt=-\int_{M_{\infty}} \omega(t)\wedge\frac{\partial\eta(t,x)}{\partial t}dt
  \end{eqnarray*}
  it the sense of currents.

We study the conical K\"ahler-Ricci flow $(\ref{TCKRF1})$ by using the following twisted K\"ahler-Ricci flow with weak initial data $\omega_{\varphi_{0}}$.
\begin{eqnarray}\label{TKRF1}
\begin{cases}
  \frac{\partial \omega_{\varepsilon}(t)}{\partial t}=-Ric(\omega_{\varepsilon}(t))+\mu\omega_{\varepsilon}(t)+\theta_{\varepsilon},\\
  \\
  \omega_{\varepsilon}(t)|_{t=0}=\omega_{\varphi_{0}},\\
  \end{cases}
  \end{eqnarray}
where $\theta_{\varepsilon}=(1-\beta)(\lambda\omega_{0}+\sqrt{-1}\partial\overline{\partial}\log(\varepsilon^{2}+|s|_{h}^{2}))$ is a smooth closed positive (1,1)-form, $s$ is the definition section of $D$ and $h$ is a smooth Hermitian metric on $-\lambda K_{M}$ with curvature $\lambda\omega_{0}$. The smooth case of the twisted K\"ahler-Ricci flow was studied in \cite{TC, GEDWA, SWFTZ, SWFTZ1, JWL, LW, JWLXZ, LMSH1, YQW}, etc.

There are some important results on the K\"ahler-Ricci flow (as well as its twisted versions with smooth twisting form) from weak initial data, such as Chen-Ding \cite{CD},  Chen-Tian \cite{CGT}, Chen-Tian-Zhang \cite{CTZ}, Guedj-Zeriahi \cite{VGAZ2},  Lott-Zhang \cite{LZ}, Nezza-Lu \cite{NL2014}, Song-Tian \cite{JSGT},  Sz\'ekelyhidi-Tosatti \cite{GSVT}, Z. Zhang \cite{ZZ}. Here, we study the K\"ahler-Ricci flow which is twisted by non-smooth twisting form, i.e. the flow (\ref{TCKRF1}).  We first obtain the long-time existence, uniqueness and regularity of the flow $(\ref{TKRF1})$ by following Song-Tian's arguments in \cite{JSGT}. Then we study the long-time existence of the conical K\"ahler-Ricci flow $(\ref{TCKRF1})$ by approximating method. In this process, in addition to getting the locally uniform regularity of the twisted K\"ahler-Ricci flow $(\ref{TKRF1})$, the most important step is to prove that $\varphi(t)$ converges to $\varphi_{0}$ in $L^{\infty}$-norm as $t\rightarrow0^{+}$  ( i.e the 4th property in Definition \ref{04.5}), where $\varphi(t)$ is a metric potential of $\omega(t)$ with respect to metric $\omega_{0}$. Here we need a new idear because of  Song-Tian's method in \cite{JSGT} is invalid. At the same time, we prove the uniqueness of the conical K\"ahler-Ricci flow by Jeffres' trick \cite{TJEF} and an improvement of the arguments in \cite{YQW}. In fact, we obtain the following theorem.

\begin{thm}\label{04} Let $M$ be a Fano manifold with complex dimension $n$, $\omega_{0}\in c_{1}(M)$ be a smooth K\"ahler metric on $M$, divisor $D\in|-\lambda K_{M}|$ $(\lambda\in\mathbb{Q})$ and $\hat{\omega}\in c_{1}(M)$ be a K\"ahler current which admits $L^{p}$ density with respect to $\omega_{0}^{n}$ for some $p>1$ and $\int_{M}\hat{\omega}^{n}=\int_{M}\omega_{0}^{n}$. Then for any $\beta\in(0,1)$, there exists a unique solution $\omega(t,\cdot)$ to the conical K\"ahler-Ricci flow $(\ref{TCKRF1})$ with weak initial data $\hat{\omega}$.
\end{thm}

Then we consider the convergence of the conical K\"ahler-Ricci flow $(\ref{TCKRF1})$. When $\lambda>0$ and there is no nontrivial holomorphic field which is tangent to $D$ along $D$, Tian-Zhu \cite{GTXHZ05} proved a Moser-Trudinger type inequality for conical K\"ahler-Einstein manifold and gave a new proof of Donaldson's openness theorem \cite{SD2}. Using the Moser-Trudinger type inequality in \cite{GTXHZ05} and following the arguments in \cite{JWLXZ}, we obtain the following convergence result of the conical K\"ahler-Ricci flow $(\ref{TCKRF1})$.
\begin{thm}\label{05}Assume that $\lambda >0$ and there is no nontrivial holomorphic field on $M$ tangent to $D$, if there exists a conical K\"ahler-Einstein metric with cone angle $2\pi\beta$ ( $0<\beta<1$) along $D$, then the conical K\"ahler-Ricci flow $(\ref{TCKRF1})$ must converge to this conical K\"ahler-Einstein metric in $C_{loc}^{\infty}$ topology outside divisor $D$ and globally in the sense of currents on $M$.
\end{thm}

\begin{rem} In this paper, we only study the convergence with positive twisted first Chern class, i.e. $\mu=1-(1-\beta)\lambda >0$. When $\mu\leq0$, one can also get the convergence of the conical K\"ahler-Ricci flow by following Chen-Wang's argument in \cite{CW1}.
\end{rem}

\medskip

The paper is organized as follows. In section $2$, we prove the long-time existence and uniqueness of the twisted K\"ahler-Ricci flow $(\ref{TKRF1})$ by adopting Song-Tian's methods in \cite{JSGT}. In section $3$,  we obtain the existence of a long-time solution to the conical K\"ahler-Ricci flow (\ref{TCKRF1}) by limiting the twisted K\"ahler-Ricci flows,  and prove that $\varphi(t)$ converges to $\varphi_{0}$ in $L^{\infty}$-norm as $t\rightarrow0^{+}$, where $\varphi(t)$ is a metric potential of $\omega(t)$ with respect to the metric $\omega_{0}$. We also prove the uniqueness of the conical K\"ahler-Ricci flow with weak initial data $\omega_{\varphi_{0}}$. In section $4$, by using the uniform Perelman's estimates along the twisted K\"ahler-Ricci flows obtained in \cite{JWLXZ}, we prove the convergence theorem under the assumptions in Theorem \ref{05}.
\medskip

{\bf  Acknowledgement:} The  authors would like to thank Professor Jiayu Li  for providing  many suggestions and encouragements. The first author also would like to thank Professor Xiaohua Zhu  for his constant help and encouragements. The second author is partially supported by the NSFC Grants 11131007 and 11571332. 

\section{The long-time existence of the twisted K\"ahler-Ricci flow with weak initial data}
\setcounter{equation}{0}

In this section, we prove the long-time existence and uniqueness of the twisted K\"ahler-Ricci flow $(\ref{TKRF1})$ by following Song-Tian's arguments in \cite{JSGT}. For further consideration in the next section, we shall pay attention to the estimates which are independent of $\varepsilon$. In the following arguments, for the sake of brevity, we only consider the flow $(\ref{TKRF1})$ in the case of $\lambda=1$  (i.e. $\mu=\beta$), where $\beta\in (0,1)$. Our arguments are also valid for any $\lambda$, only if the coefficient $\beta$ before $\omega_\varepsilon(t)$ in the case of $\lambda=1$ is replaced by $\mu=1-(1-\beta)\lambda$. We denote
 \begin{eqnarray*}
 F=\frac{(\omega_{0}+\sqrt{-1}\partial\bar{\partial}\varphi_{0})^{n}}{\omega_{0}^{n}}\in L^{p}(M,\omega_{0}^{n})\ for\ p>1.
 \end{eqnarray*}
Recall that $C^{\infty}(M)$ is dense in $L^{p}(M,\omega_{0}^{n})$. Therefore there exists a sequence of positive functions $F_{j}\in C^{\infty}(M)$ such that $\int_{M}F_{j}\omega_{0}^{n}=\int_{M}\omega_{0}^{n}$ and
\begin{eqnarray*}
\lim\limits_{j\rightarrow\infty}\|F_{j}-F\|_{L^{p}(M)}=0.
\end{eqnarray*}
By considering the complex Monge-Amp\`ere equation
\begin{eqnarray}
(\omega_{0}+\sqrt{-1}\partial\bar{\partial}\varphi_{0,j})^{n}=F_{j}\omega_{0}^{n}
\end{eqnarray}
and using the stability theorem in \cite{K00} (see also \cite{SDZZ} or \cite{VGAZ1}), we have
\begin{eqnarray}\label{00000}
\lim\limits_{j\rightarrow\infty}\|\varphi_{0,j}-\varphi_{0}\|_{L^{\infty}(M)}=0,
\end{eqnarray}
where $\varphi_{0,j}\in PSH(M,\omega_0)\cap C^\infty(M)$ satisfy
$\sup\limits_{M}(\varphi_0-\varphi_{0,j})=\sup\limits_{M}(\varphi_{0,j}-\varphi_{0})$.

Let $\omega_{\varphi_{0,j}}=\omega_0+\sqrt{-1}\partial\overline{\partial}\varphi_{0,j}$. We prove the long-time existence of the twisted K\"ahler-Ricci flow $(\ref{TKRF1})$ by using a sequence of smooth twisted K\"ahler-Ricci flows
\begin{eqnarray}\label{TKRF2}
\begin{cases}
  \frac{\partial \omega_{\varepsilon,j}(t)}{\partial t}=-Ric(\omega_{\varepsilon,j}(t))+\beta\omega_{\varepsilon,j}(t)+\theta_{\varepsilon}.\\
  \\
  \omega_{\varepsilon,j}(t)|_{t=0}=\omega_{\varphi_{0,j}}\\
  \end{cases}
\end{eqnarray}

Since the twisted K\"ahler-Ricci flow preserves the K\"ahler class, we can write the flow $(\ref{TKRF2})$ as the parabolic Monge-Amp\'ere equation on potentials,
\begin{eqnarray}\label{CMAE1}
\begin{cases}
  \frac{\partial \varphi_{\varepsilon,j}(t)}{\partial t}=\log\frac{(\omega_{0}+\sqrt{-1}\partial\bar{\partial}\varphi_{\varepsilon,j}(t))^{n}}{\omega_{0}^{n}}+F_{0}\\
 \ \ \ \ \ \ \ \ \ \ \ \ \ +\beta\varphi_{\varepsilon,j}(t)+\log(\varepsilon^{2}+|s|_{h}^{2})^{1-\beta},\\
  \varphi_{\varepsilon,j}(0)=\varphi_{0,j}\\
  \end{cases}
\end{eqnarray}
where $F_{0}$ satisfies $-Ric(\omega_{0})+\omega_{0}=\sqrt{-1}\partial\overline{\partial}F_{0}$, $\frac{1}{V}\int_{M}e^{-F_{0}}dV_{0}=1$ and $dV_{0}=\frac{\omega_{0}^{n}}{n!}$. By using the function
\begin{eqnarray}\chi(\varepsilon^{2}+|s|_h^2)=\frac{1}{\beta}\int_{0}^{|s|_h^2}\frac{(\varepsilon^{2}+r)^{\beta}-
\varepsilon^{2\beta}}{r}dr
\end{eqnarray}
which was given by Campana-Guenancia-P$\breve{a}$un in \cite{CGP}, we can rewrite the flow (\ref{CMAE1}) as
\begin{eqnarray}\label{CMAE2}
\begin{cases}
  \frac{\partial \phi_{\varepsilon,j}(t)}{\partial t}=\log\frac{(\omega_{\varepsilon}+\sqrt{-1}\partial\bar{\partial}\phi_{\varepsilon,j}(t))^{n}}{\omega_{\varepsilon}^{n}}+F_{\varepsilon}+
  \beta(\phi_{\varepsilon,j}(t)+k\chi(\varepsilon^2+|s|_h^2)),\\
  \\
  \phi_{\varepsilon,j}(0)=\varphi_{0,j}-k\chi(\varepsilon^2+|s|_h^2):=\phi_{\varepsilon,0,j}\\
  \end{cases}
\end{eqnarray}
where $\phi_{\varepsilon,j}(t)=\varphi_{\varepsilon,j}(t)-k\chi(\varepsilon^2+|s|_h^2)$, $\omega_{\varepsilon}=\omega_{0}+\sqrt{-1}k\partial\overline{\partial}\chi(\varepsilon^{2}+|s|_{h}^{2})$, $F_{\varepsilon}=F_{0}+\log(\frac{\omega_{\varepsilon}^{n}}{\omega_{0}^{n}}\cdot(\varepsilon^{2}+|s|_{h}^{2})^{1-\beta})$. We know that $\chi(\varepsilon^2+|s|_h^2)$ and $F_{\varepsilon}$ are uniformly bounded (see $(15)$ and $(25)$ in \cite{CGP}).

\begin{pro}\label{200} For any $T>0$, there exists constant $C$ depending only on $\|\varphi_0\|_{L^\infty(M)}$, $\beta$, $n$, $\omega_{0}$ and $T$ such that for any $t\in[0,T]$, $\varepsilon>0$ and $j\in\mathbb{N}^+$,
\begin{eqnarray}\|\phi_{\varepsilon,j}(t)\|_{L^\infty(M)}\leq C.
\end{eqnarray}
Furthermore, for any $j, l$, we have
\begin{eqnarray}\label{201}
\|\phi_{\varepsilon,j}(t)-\phi_{\varepsilon,l}(t)\|_{L^\infty([0,T]\times M)}\leq e^{\beta T}\|\varphi_{0,j}-\varphi_{0,l}\|_{L^\infty(M)}.
\end{eqnarray}
In particular, $\{\varphi_{\varepsilon,j}(t)\}$ satisfies
\begin{eqnarray}\label{2001}\lim\limits_{j.l\rightarrow\infty}\|\varphi_{\varepsilon,j}(t)-\varphi_{\varepsilon,l}(t)\|_{L^{\infty}([0,T]\times M)}=0.
\end{eqnarray}
\end{pro}

{\bf Proof:}\ \ From equation $(\ref{CMAE2})$, we have
\begin{eqnarray*}
\frac{\partial e^{-\beta t}\phi_{\varepsilon,j}(t)}{\partial t}&=&e^{-\beta t}\log\frac{(e^{-\beta t}\omega_{\varepsilon}+\sqrt{-1}\partial\bar{\partial}e^{-\beta t}\phi_{\varepsilon,j}(t))^{n}}{(e^{-\beta t}\omega_{\varepsilon})^{n}}\\
&\ &+e^{-\beta t}(F_{\varepsilon}+k\beta\chi(\varepsilon^2+|s|_h^2))\\
&\leq&e^{-\beta t}\log\frac{(e^{-\beta t}\omega_{\varepsilon}+\sqrt{-1}\partial\bar{\partial}e^{-\beta t}\phi_{\varepsilon,j}(t))^{n}}{(e^{-\beta t}\omega_{\varepsilon})^{n}}+Ce^{-\beta t},
\end{eqnarray*}
which is equivalent to
\begin{eqnarray*}
\frac{\partial }{\partial t}\big(e^{-\beta t}(\phi_{\varepsilon,j}(t)+\frac{C}{\beta})\big)\leq e^{-\beta t}\log\frac{\big(e^{-\beta t}\omega_{\varepsilon}+\sqrt{-1}\partial\bar{\partial}e^{-\beta t}(\phi_{\varepsilon,j}(t)+\frac{C}{\beta})\big)^{n}}{(e^{-\beta t}\omega_{\varepsilon})^{n}},
\end{eqnarray*}
where constant $C$ depends only on $\|\varphi_0\|_{L^\infty(M)}$, $\beta$, $n$ and $\omega_{0}$.

For any $\delta>0$, we denote $\tilde{\phi}_{\varepsilon,j}(t)=e^{-\beta t}(\phi_{\varepsilon,j}(t)+\frac{C}{\beta})-\delta t$. Let $(t_0,x_0)$ be the maximum point of $\tilde{\phi}_{\varepsilon,j}(t)$ on $[0,T]\times M$. If $t_0>0$, by maximum principle, we have
\begin{eqnarray*}
0&\leq&\frac{\partial }{\partial t}\big(e^{-\beta t}(\phi_{\varepsilon,j}(t)+\frac{C}{\beta})\big)(t_0,x_0)-\delta\\
&\leq& e^{-\beta t}\log\frac{\big(e^{-\beta t}\omega_{\varepsilon}+\sqrt{-1}\partial\bar{\partial}\tilde{\phi}_{\varepsilon,j}(t)\big)^{n}}{(e^{-\beta t}\omega_{\varepsilon})^{n}}(t_0,x_0)-\delta\\
&\leq&-\delta,
\end{eqnarray*}
which is impossible. Hence $t_0=0$, then
\begin{eqnarray*}
\phi_{\varepsilon,j}(t)\leq e^{\beta t}\sup\limits_{M}\phi_{\varepsilon,j}(0)+\delta Te^{\beta T}+\frac{C}{\beta}(e^{\beta T}-1).
\end{eqnarray*}
Let $\delta\rightarrow0$, we obtain
\begin{eqnarray}\label{202}
\phi_{\varepsilon,j}(t)\leq e^{\beta t}\sup\limits_{M}\phi_{\varepsilon,j}(0)+\frac{C}{\beta}(e^{\beta T}-1).
\end{eqnarray}
By the same arguments, we can get the lower bound of $\phi_{\varepsilon,j}(t)$
\begin{eqnarray}\label{203}
\phi_{\varepsilon,j}(t)\geq e^{\beta t}\inf\limits_{M}\phi_{\varepsilon,j}(0)-\frac{C}{\beta}(e^{\beta T}-1).
\end{eqnarray}
Combining $(\ref{202})$ and $(\ref{203})$, we have
\begin{eqnarray*}
\|\phi_{\varepsilon,j}(t)\|_{L^\infty(M)}\leq e^{\beta T}\|\phi_{\varepsilon,j}(0)\|_{L^\infty(M)}+\frac{C}{\beta}(e^{\beta T}-1)\leq C,
\end{eqnarray*}
where constant $C$ depends only on $\|\varphi_0\|_{L^\infty(M)}$, $\beta$, $n$, $\omega_{0}$ and $T$.

Let $\psi_{\varepsilon,j,l}(t)=\phi_{\varepsilon,j}(t)-\phi_{\varepsilon,l}(t)$, then $\psi_{\varepsilon,j,l}$ satisfies the following equation
\begin{eqnarray}
\begin{cases}
  \frac{\partial \psi_{\varepsilon,j,l}(t)}{\partial t}=\log\frac{\big(\omega_{\varepsilon}+\sqrt{-1}\partial\bar{\partial}\phi_{\varepsilon,l}(t)+\sqrt{-1}\partial\bar{\partial}\psi_{\varepsilon,j,l}(t)\big)^{n}}
  {(\omega_{\varepsilon}+\sqrt{-1}\partial\bar{\partial}\phi_{\varepsilon,l}(t))^{n}}+
  \beta\psi_{\varepsilon,j,l}(t).\\
  \\
  \psi_{\varepsilon,j,l}(0)=\varphi_{0,j}-\varphi_{0,l}\\
  \end{cases}
\end{eqnarray}
By the same arguments as that in the first part, we have
\begin{eqnarray*}
\|\psi_{\varepsilon,j,l}(t)\|_{L^\infty([0,T]\times M)}\leq e^{\beta T}\|\varphi_{0,j}-\varphi_{0,l}\|_{L^\infty(M)}.
\end{eqnarray*}
Since $\{\varphi_{0,j}\}$ is a Cauchy sequence in $L^\infty$-norm, we conclude the limit $(\ref{2001})$.\QEDB

\medskip

We now prove the uniform equivalence of the volume forms along the complex Monge-Amp\`ere flow $(\ref{CMAE2})$.

\begin{lem}\label{204} For any $T>0$, there exists constant $C$ depending only on $\|\varphi_0\|_{L^\infty(M)}$, $n$, $\beta$, $\omega_{0}$ and $T$ such that for any $t\in(0,T]$, $\varepsilon>0$ and $j\in\mathbb{N}^+$,
\begin{eqnarray}\frac{t^n}{C}\leq\frac{(\omega_{\varepsilon}+\sqrt{-1}\partial\bar{\partial}\phi_{\varepsilon,j}(t))^{n}}{\omega_{\varepsilon}^n}\leq e^{\frac{C}{t}}.
\end{eqnarray}
\end{lem}

\medskip

{\bf Proof:}\ \ Let $\Delta_{\varepsilon,j}$ be the Laplacian operator associated to the K\"ahler form $\omega_{\varepsilon,j}(t)=\omega_{\varepsilon}+\sqrt{-1}\partial\bar{\partial}\phi_{\varepsilon,j}(t)$. Straightforward calculations show that
\begin{eqnarray}(\frac{\partial}{\partial t}-\Delta_{\varepsilon,j})\dot{\phi}_{\varepsilon,j}(t)=\beta\dot{\phi}_{\varepsilon,j}(t).
\end{eqnarray}

Let $H_{\varepsilon,j}^+(t)=t\dot{\phi}_{\varepsilon,j}(t)-A\phi_{\varepsilon,j}(t)$, where $A$ is a sufficiently large number (for example $A=\beta T+2$). Then $H_{\varepsilon,j}^+(0)=-A\phi_{\varepsilon,j}(0)$ is uniformly bounded by a constant $C$ depending only on $\|\varphi_0\|_{L^\infty(M)}$, $\beta$, $n$, $\omega_{0}$ and $T$.
\begin{eqnarray}\nonumber(\frac{\partial}{\partial t}-\Delta_{\varepsilon,j})H_{\varepsilon,j}^+(t)&=&(1+\beta t-A)\dot{\phi}_{\varepsilon,j}(t)+A\Delta_{\varepsilon,j}\phi_{\varepsilon,j}(t)\\
&\leqslant&(1+\beta t-A)\dot{\phi}_{\varepsilon,j}(t)+An.
\end{eqnarray}
By the maximum principle, $H_{\varepsilon,j}^+(t)$ is uniformly bounded from above by a constant $C$ depending only on $\|\varphi_0\|_{L^\infty(M)}$, $n$, $\beta$, $\omega_{0}$ and $T$.

Let $H_{\varepsilon,j}^-(t)=\dot{\phi}_{\varepsilon,j}(t)+\phi_{\varepsilon,j}(t)-n\log t$. Then $H^-(t)$ tends to $+\infty$ as $t\rightarrow0^+$ and
\begin{eqnarray}(\frac{\partial}{\partial t}-\Delta_{\varepsilon,j})H_{\varepsilon,j}^-(t)=(\beta+1)\dot{\phi}_{\varepsilon,j}(t)+tr_{\omega_{\varepsilon,j}(t)}\omega_\varepsilon-n-\frac{n}{t}.
\end{eqnarray}

Assume that $(t_0,x_0)$ is the minimum point of $H_{\varepsilon,j}^-(t)$ on $[0,T]\times M$. We conclude that $t_0>0$ and there exists constant $C_1$, $C_2$ and $C_3$ such that
\begin{eqnarray}\label{2002}\nonumber(\frac{\partial}{\partial t}-\Delta_{\varepsilon,j})H_{\varepsilon,j}^-(t)|_{(t_0,x_0)}&\geq& \big(C_1(\frac{\omega_{\varepsilon}^n}{\omega^{n}_{\varepsilon,j}(t)})^{\frac{1}{n}}+C_2\log\frac{\omega^n_{\varepsilon,j}(t)}{\omega_{\varepsilon}^n}
-\frac{C_3}{t}\big)|_{(t_0,x_0)}\\
&\geq&\big(\frac{C_1}{2}(\frac{\omega_{\varepsilon}^n}{\omega^{n}_{\varepsilon,j}(t)})^{\frac{1}{n}}
-\frac{C_3}{t}\big)|_{(t_0,x_0)},
\end{eqnarray}
where constant $C_1$ depends only on $n$, $C_2$ depends only on $\beta$ and $C_3$ depends only on $n$, $\omega_{0}$, $\|\varphi_0\|_{L^\infty(M)}$, $\beta$ and $T$. In inequality $(\ref{2002})$, without loss of generality, we assume that $\frac{\omega_{\varepsilon}^n}{\omega^{n}_{\varepsilon,j}(t)}>1$ and $\frac{C_1}{2}(\frac{\omega_{\varepsilon}^n}{\omega^{n}_{\varepsilon,j}(t)})^{\frac{1}{n}}+C_2\log\frac{\omega^n_{\varepsilon,j}(t)}
{\omega_{\varepsilon}^n}\geq0$ at $(t_0,x_0)$. By the maximum principle, we have
\begin{eqnarray}\omega^{n}_{\varepsilon,j}(t_0,x_0)\geq C_4t^n \omega_{\varepsilon}^n(x_0),
\end{eqnarray}
where $C_4$ independent of $\varepsilon$ and $j$. Then it easily follows that $H_{\varepsilon,j}^-(t)$ is bounded from below by a constant $C$ depending only on $\|\varphi_0\|_{L^\infty(M)}$, $n$, $\beta$, $\omega_{0}$ and $T$.\QEDB

\medskip

In the following lemma, we prove the uniform equivalence of the metrics along the twisted K\"ahler-Ricci flow $(\ref{TKRF2})$.

\begin{lem}\label{205} For any $T>0$, there exists constant $C$ depending only on $\|\varphi_0\|_{L^\infty(M)}$, $n$, $\beta$, $\omega_{0}$ and $T$ such that for any $t\in(0,T]$, $\varepsilon>0$ and $j\in\mathbb{N}^+$,
\begin{eqnarray}e^{-\frac{C}{t}}\omega_\varepsilon\leq\omega_{\varepsilon,j}(t)\leq e^{\frac{C}{t}}\omega_\varepsilon.
\end{eqnarray}
\end{lem}

{\bf Proof:}\ \ Let
\begin{eqnarray}\Psi_{\varepsilon,\rho}=B\frac{1}{\rho}\int_{0}^{|s|_h^2}\frac{(\varepsilon^{2}+r)^{\rho}-
\varepsilon^{2\rho}}{r}dr
\end{eqnarray}
be the uniform bound function introduced by Guenancia-P$\breve{a}$un in \cite{GP1}. By choosing suitable $B$ and $\rho$, and following the arguments in section $2$ of \cite{JWLXZ}, we have
\begin{eqnarray}\nonumber&\ &(\frac{\partial}{\partial t}-\Delta_{\varepsilon,j})(t\log tr_{\omega_{\varepsilon}}\omega_{\varepsilon,j}(t)+t\Psi_{\varepsilon,\rho})\\
&\leq&\log tr_{\omega_{\varepsilon}}\omega_{\varepsilon,j}(t)+Ctr_{\omega_{\varepsilon,j}(t)}\omega_{\varepsilon}+C,
\end{eqnarray}
where constant $C$ depends only on $n$, $\beta$, $\omega_{0}$ and $T$.

Let $H_{\varepsilon,j}(t)=t\log tr_{\omega_{\varepsilon}}\omega_{\varepsilon,j}(t)+t\Psi_{\varepsilon,\rho}-A\phi_{\varepsilon,j}(t)$, $A$ be a sufficiently large constant and $(t_0,x_0)$ be the maximum point of $H_{\varepsilon,j}(t)$ on $[0,T]\times M$. We need only consider $t_0>0$. By the inequality
\begin{eqnarray}\label{206}
tr_{\omega_{\varepsilon}}\omega_{\varepsilon,j}(t)\leq\frac{1}{(n-1)!}(tr_{\omega_{\varepsilon,j}(t)}\omega_{\varepsilon})^{n-1}\frac{\omega^n_{\varepsilon,j}(t)}{\omega_{\varepsilon}^n},
\end{eqnarray}
we conclude that
\begin{eqnarray}\nonumber(\frac{\partial}{\partial t}-\Delta_{\varepsilon,j})H_{\varepsilon,j}(t)&\leq&\log tr_{\omega_{\varepsilon}}\omega_{\varepsilon,j}(t)+Ctr_{\omega_{\varepsilon,j}(t)}\omega_{\varepsilon}-A\dot{\phi}_{\varepsilon,j}(t)-Atr_{\omega_{\varepsilon,j}(t)}\omega_{\varepsilon}+C\\\nonumber
&\leq&(n-1)\log tr_{\omega_{\varepsilon,j}(t)}\omega_{\varepsilon}-\frac{A}{2}tr_{\omega_{\varepsilon,j}(t)}\omega_{\varepsilon}-(A-1)\log\frac{\omega^n_{\varepsilon,j}(t)}{\omega_{\varepsilon}^n}+C,
\end{eqnarray}
where constant $C$ depends only on $\|\varphi_0\|_{L^\infty(M)}$, $n$, $\beta$, $\omega_{0}$ and $T$.

Without loss of generality, we assume that
$-\frac{A}{4}tr_{\omega_{\varepsilon,j}(t)}\omega_{\varepsilon}+(n-1)\log tr_{\omega_{\varepsilon,j}(t)}\omega_{\varepsilon}\leq0$ at $(t_0,x_0)$. Then at $(t_0,x_0)$, by Lemma \ref{204}, we have
\begin{eqnarray}(\frac{\partial}{\partial t}-\Delta_{\varepsilon,j})H_{\varepsilon,j}(t)\leq-\frac{A}{4}tr_{\omega_{\varepsilon,j}(t)}\omega_{\varepsilon}-C\log t+C.
\end{eqnarray}
By the maximum principle, at $(t_0,x_0)$,
\begin{eqnarray}tr_{\omega_{\varepsilon,j}(t)}\omega_{\varepsilon}\leq C\log \frac{1}{t}+C.
\end{eqnarray}
By using inequality $(\ref{206})$, at $(t_0,x_0)$,
\begin{eqnarray}tr_{\omega_{\varepsilon}}\omega_{\varepsilon,j}(t)\leq C(\log \frac{1}{t}+1)^{n-1}e^{\frac{C}{t}}\leq e^{\frac{2C}{t}},
\end{eqnarray}
where constant $C$ depends only on $\|\varphi_0\|_{L^\infty(M)}$, $n$, $\beta$, $\omega_{0}$ and $T$. Hence we have
\begin{eqnarray}\label{207}tr_{\omega_{\varepsilon}}\omega_{\varepsilon,j}(t)\leq e^{\frac{C}{t}}
\end{eqnarray}
for some constant $C$ depending only on $\|\varphi_0\|_{L^\infty(M)}$, $n$, $\beta$, $\omega_{0}$ and $T$.

Furthermore, by inequality $(\ref{206})$ again, we know
\begin{eqnarray}\label{208}tr_{\omega_{\varepsilon,j}(t)}\omega_{\varepsilon}\leq e^{\frac{C}{t}},
\end{eqnarray}
where constant $C$ depends only on $\|\varphi_0\|_{L^\infty(M)}$, $n$, $\beta$, $\omega_{0}$ and $T$. From $(\ref{207})$ and $(\ref{208})$, we prove the lemma.\QEDB

By Lemma \ref{205} and the fact that $\omega_\varepsilon>\gamma \omega_0$ for some uniform constant $\gamma$ (see inequality $(24)$ in \cite{CGP}), we have
\begin{eqnarray}e^{-\frac{C}{t}}\omega_0\leq\omega_{\varepsilon,j}(t)\leq C_\varepsilon e^{\frac{C}{t}}\omega_0,
\end{eqnarray}
on $(0,T]\times M$, where $C$ is a uniform constant and $C_\varepsilon$ depends on $\varepsilon$. We next prove the Calabi's $C^3$-estimates. Denote
\begin{eqnarray}S_{\varepsilon,j}=|\nabla_{\omega_0}\omega_{\varepsilon,j}(t)|^{2}_{\omega_{\varepsilon,j}(t)}=
g_{\varepsilon,j}^{i\bar{m}}g_{\varepsilon,j}^{k\bar{l}}g_{\varepsilon,j}^{p\bar{q}}\nabla_{0 i}(g_{\varepsilon,j}) _{k\bar{q}}\overline{\nabla}_{0 m}(g_{\varepsilon,j})_{ p\bar{l}}.
\end{eqnarray}

\begin{lem}\label{209} For any $T>0$ and $\varepsilon>0$, there exist constants $C_\varepsilon$ and $C$ such that for any $t\in(0,T]$ and $j\in\mathbb{N}^+$,
\begin{eqnarray}S_{\varepsilon,j}\leq C_\varepsilon e^{\frac{C}{t}},
\end{eqnarray}
where constant $C$ depends only on $\|\varphi_0\|_{L^\infty(M)}$, $n$, $\beta$, $\omega_{0}$ and $T$, and constant $C_\varepsilon$ depends in addition on $\varepsilon$.
\end{lem}

\medskip

{\bf Proof:}\ \ By the similar arguments in \cite{LW} or \cite{JWLXZ} and choosing sufficiently large $\alpha$ and $\beta$, we have
\begin{eqnarray}(\frac{\partial}{\partial t}-\Delta_{\varepsilon,j})(e^{-\frac{2\alpha}{t}}S_{\varepsilon,j})&\leq& C_\varepsilon e^{-\frac{\alpha}{t}}S_{\varepsilon,j}+C_\varepsilon,\\
(\frac{\partial}{\partial t}-\Delta_{\varepsilon,j})(e^{-\frac{2\gamma}{t}}tr_{\omega_0}\omega_{\varepsilon,j}(t))&\leq& C_\varepsilon -C_\varepsilon^{-1}e^{-\frac{3\gamma}{t}}S_{\varepsilon,j}.
\end{eqnarray}
By choosing $A_\varepsilon=C_\varepsilon(C_\varepsilon+1)$ and $\alpha=3\gamma$,
\begin{eqnarray}(\frac{\partial}{\partial t}-\Delta_{\varepsilon,j})(e^{-\frac{2\alpha}{t}}S_{\varepsilon,j}+A_\varepsilon e^{-\frac{2\gamma}{t}}tr_{\omega_0}\omega_{\varepsilon,j}(t))\leq - e^{-\frac{3\gamma}{t}}S_{\varepsilon,j}+C_\varepsilon.
\end{eqnarray}
By the maximum principle, we have
\begin{eqnarray}S_{\varepsilon,j}\leq C_\varepsilon e^{\frac{C}{t}}\ \ \ \ on\ \ \ (0,T]\times M
\end{eqnarray}
for some constant $C$ depending only on $\|\varphi_0\|_{L^\infty(M)}$, $n$, $\beta$, $\omega_{0}$ and $T$, and constant $C_\varepsilon$ depending in addition on $\varepsilon$.\QEDB

\medskip

By using the Schauder regularity theory and equation $(\ref{CMAE1})$, we get the high order estimates of $\varphi_{\varepsilon,j}(t)$.

\begin{pro}\label{210} For any $0<\delta<T<\infty$, $\varepsilon>0$ and $k\geq0$, there exists constant $C_{\varepsilon,\delta,T,k}$ depending only on $\delta$, $T$, $\varepsilon$, $k$, $n$, $\beta$, $\omega_{0}$ and $\|\varphi_0\|_{L^\infty(M)}$, such that for any $j\in\mathbb{N}^+$,
\begin{eqnarray}\|\varphi_{\varepsilon,j}(t)\|_{C^{k}\big([\delta,T]\times M\big)}\leq C_{\varepsilon,\delta,T,k}.
\end{eqnarray}
\end{pro}

By $(\ref{2001})$, for any $T>0$, $\varphi_{\varepsilon,j}(t)$ converges to $\varphi_\varepsilon(t)\in L^\infty([0,T]\times M)$ uniformly in $L^\infty([0,T]\times M)$. For any $0<\delta<T<\infty$ and $\varepsilon>0$, $\varphi_{\varepsilon,j}(t)$ is uniformly bounded (depends on $\varepsilon$) in $C^{\infty}([\delta,T]\times M)$. Therefore $\varphi_{\varepsilon,j}(t)$ converges to $\varphi_\varepsilon(t)$ in $C^\infty([\delta,T]\times M)$. Hence for any $\varepsilon>0$, $\varphi_\varepsilon(t)\in C^\infty((0,\infty)\times M)$.

\begin{pro}\label{211} For any $\varepsilon>0$, $\varphi_\varepsilon(t)\in C^0([0,\infty)\times M)$ and
\begin{eqnarray}\label{212}\lim\limits_{t\rightarrow0^+}\|\varphi_{\varepsilon}(t)-\varphi_{0}\|_{L^{\infty}(M)}=0.
\end{eqnarray}
\end{pro}

{\bf Proof:}\ \ For any $(t,z)\in(0,T]\times M$,
\begin{eqnarray}\nonumber|\varphi_{\varepsilon}(t,z)-\varphi_{0}(z)|&\leq&|\varphi_{\varepsilon}(t,z)-\varphi_{\varepsilon,j}(t,z)|
+|\varphi_{\varepsilon,j}(t,z)-\varphi_{0,j}(z)|\\
&\ &+|\varphi_{0,j}(z)-\varphi_{0}(z)|.
\end{eqnarray}
Since $\varphi_{\varepsilon,j}(t)$ is a Cauchy sequence in $L^\infty([0,T]\times M)$,
\begin{eqnarray}\lim\limits_{j\rightarrow\infty}\|\varphi_{\varepsilon}(t,z)-\varphi_{\varepsilon,j}(t,z)\|_{L^{\infty}([0,T]\times M)}=0.
\end{eqnarray}
From $(\ref{00000})$, we have
\begin{eqnarray}
\lim\limits_{j\rightarrow\infty}\|\varphi_{0,j}(z)-\varphi_{0}(z)\|_{L^{\infty}(M)}=0,
\end{eqnarray}
For any $\epsilon>0$, there exists $N$ such that for any $j>N$,
\begin{eqnarray*}\sup\limits_{[0,T]\times M}|\varphi_{\varepsilon}(t,z)-\varphi_{\varepsilon,j}(t,z)|&<&\frac{\epsilon}{3},\\
\sup\limits_{M}|\varphi_{0,j}(z)-\varphi_{0}(z)|&<&\frac{\epsilon}{3}.
\end{eqnarray*}
On the other hand, fix such $j$, there exists $0<\delta<T$ such that
\begin{eqnarray}\sup\limits_{[0,\delta]\times M}|\varphi_{\varepsilon,j}(t,z)-\varphi_{0,j}|<\frac{\epsilon}{3}.
\end{eqnarray}
Combining the above estimates together, for any $t\in[0,\delta]$ and $z\in M$,
\begin{eqnarray}|\varphi_{\varepsilon}(t,z)-\varphi_{0}(z)|<\epsilon.
\end{eqnarray}
This completes the proof of the lemma.\QEDB

\begin{pro} \label{214}$\varphi_\varepsilon(t)$ is the unique solution to the parabolic Monge-Amp\`ere equation
\begin{eqnarray}\label{CMAE3}
\begin{cases}
  \frac{\partial \varphi_{\varepsilon}(t)}{\partial t}=\log\frac{(\omega_{0}+\sqrt{-1}\partial\bar{\partial}\varphi_{\varepsilon}(t))^{n}}{\omega_{0}^{n}}+F_{0}\\
  \ \ \ \ \ \ \ \ \ \ +\beta\varphi_{\varepsilon}(t)
  +\log(\varepsilon^{2}+|s|_{h}^{2})^{1-\beta},\ \ \ \ (0,\infty)\times M\\
  \varphi_{\varepsilon}(0)=\varphi_{0}\\
  \end{cases}
\end{eqnarray}
in the space of $C^0\big([0,\infty)\times M\big)\cap C^\infty\big((0,\infty)\times M\big)$.
\end{pro}

\medskip

{\bf Proof:}\ \ By proposition \ref{211}, we only need to prove the uniqueness. Suppose there exists another solution $\tilde{\varphi}_\varepsilon(t)\in C^0\big([0,\infty)\times M\big)\cap C^\infty\big((0,\infty)\times M\big)$ to the Monge-Amp\`ere equation $(\ref{CMAE3})$.

Let $\psi_\varepsilon(t)=\tilde{\varphi}_\varepsilon(t)-\varphi_\varepsilon(t)$. Then
\begin{eqnarray}
\begin{cases}
  \frac{\partial \psi_{\varepsilon}(t)}{\partial t}=\log\frac{\big(\omega_{0}+\sqrt{-1}\partial\bar{\partial}\varphi_{\varepsilon}(t)+\sqrt{-1}\partial\bar{\partial}\psi_{\varepsilon}(t)\big)^{n}}
  {(\omega_{0}+\sqrt{-1}\partial\bar{\partial}\varphi_{\varepsilon}(t))^{n}}+
  \beta\psi_{\varepsilon}(t).\\
  \\
  \psi_{\varepsilon}(0)=0\\
  \end{cases}
\end{eqnarray}
For any $T>0$, by the same arguments as that in the proof of Proposition \ref{200}, we have
\begin{eqnarray*}
\|\psi_{\varepsilon}(t)\|_{L^\infty([0,T]\times M)}\leq e^{\beta T}\|\psi_{\varepsilon}(0)\|_{L^\infty(M)}=0.
\end{eqnarray*}
Hence $\psi_{\varepsilon}(t)=0$, that is $\tilde{\varphi}_\varepsilon(t)=\varphi_\varepsilon(t)$.\QEDB

\medskip

By the similar arguments as that in \cite{JSGT}, we prove the uniqueness theorems of the twisted K\"ahler-Ricci flow.

\begin{thm}  Let $M$ be a Fano manifold with complex dimension $n$, $\omega_{0}\in c_{1}(M)$ be a smooth K\"ahler metric on $M$ and $\hat{\omega}\in c_{1}(M)$ be a K\"ahler current which admits $L^{p}$ density with respect to $\omega_{0}^{n}$ for some $p>1$ and $\int_{M}\hat{\omega}^{n}=\int_{M}\omega_{0}^{n}$. Then there exists a unique solution $\omega_\varepsilon(t)\in C^\infty\big((0,\infty)\times M\big)$ to the twisted K\"ahler-Ricci flow $(\ref{TKRF1})$ with initial data $\hat{\omega}$ in the following sense.

\medskip

$(1)$ $\frac{\partial \omega_{\varepsilon}(t)}{\partial t}=-Ric(\omega_{\varepsilon}(t))+\beta\omega_{\varepsilon}(t)+\theta_{\varepsilon}$ on $(0,\infty)\times M$;

\medskip

$(2)$   There exists $\varphi_\varepsilon(t)\in C^0\big([0,\infty)\times M\big)\cap C^\infty\big((0,\infty)\times M\big)$ such that $\omega_\varepsilon(t)=\omega_{0}+\sqrt{-1}\partial\bar{\partial}\varphi_{\varepsilon}(t)$ and
\begin{eqnarray}\lim\limits_{t\rightarrow0^+}\|\varphi_{\varepsilon}(t)-\varphi_{0}\|_{L^{\infty}(M)}=0,
\end{eqnarray}
where $\varphi_{0}\in\mathcal{E}_{p}(M,\omega_{0})$ is a metric potential of $\hat{\omega}$ with respect to $\omega_{0}$. In particular, $\omega_\varepsilon(t)$ converges in the sense of distribution to $\hat{\omega}$ as $t\rightarrow0$.
\end{thm}

{\bf Proof:}\ \ From Proposition \ref{214}, we know that there exists a solution $\omega_{\varepsilon}(t)=\omega_0+\sqrt{-1}\partial\bar{\partial}\varphi_{\varepsilon}(t)$ to the twisted K\"ahler-Ricci flow $(\ref{TKRF1})$ with initial data $\hat{\omega}$, where $\varphi_{\varepsilon}(t)\in C^0\big([0,\infty)\times M\big)\bigcap C^\infty\big((0,\infty)\times M\big)$ satisfies
\begin{eqnarray}\lim\limits_{t\rightarrow0^+}\|\varphi_{\varepsilon}(t)-\varphi_{0}\|_{L^{\infty}(M)}=0
\end{eqnarray}
for some metric potential $\varphi_{0}\in\mathcal{E}_{p}(M,\omega_{0})$ of $\hat{\omega}$ with respect to $\omega_{0}$. Suppose that there is another solution $\tilde{\omega}_{\varepsilon}(t)=\omega_0+\sqrt{-1}\partial\bar{\partial}\tilde{\varphi}_{\varepsilon}(t)$ to the twisted K\"ahler-Ricci flow $(\ref{TKRF1})$ with initial data $\hat{\omega}$. Then $\tilde{\varphi}_{\varepsilon}(t)\in C^0\big([0,\infty)\times M\big)\bigcap C^\infty\big((0,\infty)\times M\big)$ satisfies
\begin{eqnarray}\label{215}\nonumber\frac{\partial \tilde{\varphi}_{\varepsilon}(t)}{\partial t}&=&\log\frac{(\omega_{0}+\sqrt{-1}\partial\bar{\partial}\tilde{\varphi}_{\varepsilon}(t))^{n}}{\omega_{0}^{n}}+F_{0}\\
 &\ &+\beta\tilde{\varphi}_{\varepsilon}(t)+\log(\varepsilon^{2}+|s|_{h}^{2})^{1-\beta}+f_\varepsilon(t)
\end{eqnarray}
on $(0,\infty)\times M$ for a smooth function $f_{\varepsilon}(t)$ on $(0,\infty)$ and
\begin{eqnarray*}\lim\limits_{t\rightarrow0^{+}}\|\tilde{\varphi}_{\varepsilon}(t)-\tilde{\varphi}_{0}\|_{L^{\infty}(M)}=0,
\end{eqnarray*}
where $\tilde{\varphi}_{0}\in\mathcal{E}_{p}(M,\omega_{0})$ is also a metric potential of $\hat{\omega}$ with respect to $\omega_{0}$. At the same time, we have $\varphi_{0}=\tilde{\varphi}_{0}+\tilde{C}$.

Let $\hat{\varphi}(t)=\tilde{\varphi}(t)+\tilde{C}e^{\beta t}$. It is obvious that $\hat{\varphi}_{\varepsilon}(t)\in C^0\big([0,\infty)\times M\big)\bigcap C^\infty\big((0,\infty)\times M\big)$ is a solution to equation $(\ref{215})$ and satisfies
\begin{eqnarray*}\lim\limits_{t\rightarrow0^{+}}\|\hat{\varphi}_{\varepsilon}(t)-\varphi_{0}\|_{L^{\infty}(M)}=0.
\end{eqnarray*}

Now we consider the function $\psi_\varepsilon(t)=\hat{\varphi}_\varepsilon(t)-\varphi_\varepsilon(t)$.
\begin{eqnarray}\label{216}
\begin{cases}
  \frac{\partial \psi_{\varepsilon}(t)}{\partial t}=\log\frac{\big(\omega_{0}+\sqrt{-1}\partial\bar{\partial}\varphi_{\varepsilon}(t)+\sqrt{-1}\partial\bar{\partial}\psi_{\varepsilon}(t)\big)^{n}}
  {(\omega_{0}+\sqrt{-1}\partial\bar{\partial}\varphi_{\varepsilon}(t))^{n}}+\beta\psi_{\varepsilon}(t)+f_\varepsilon(t).\\
  \\
  \psi_{\varepsilon}(0)=0\\
  \end{cases}
\end{eqnarray}
For any $0<t_1<t_2<\infty$, by the same arguments as that in the proof of Proposition \ref{200}, we have
\begin{eqnarray*}
\sup\limits_{M}\psi_\varepsilon(t_2)&\leq&e^{\beta(t_2-t_1)}\sup\limits_{M}\psi_\varepsilon(t_1)+\int_{t_1}^{t_2}e^{\beta(t_2-t)}f_\varepsilon(t)dt,\\
\inf\limits_{M}\psi_\varepsilon(t_2)&\geq&e^{\beta(t_2-t_1)}\inf\limits_{M}\psi_\varepsilon(t_1)+\int_{t_1}^{t_2}e^{\beta(t_2-t)}f_\varepsilon(t)dt.
\end{eqnarray*}
Therefore, we obtain
\begin{eqnarray*}
\inf\limits_{M}\psi_\varepsilon(t_2)&\geq&\sup\limits_{M}\psi_\varepsilon(t_2)-e^{\beta(t_2-t_1)}(\sup\limits_{M}\psi_\varepsilon(t_1)-\inf\limits_{M}\psi_\varepsilon(t_1)).
\end{eqnarray*}
Let $t_1\rightarrow0^+$, we have
\begin{eqnarray*}
\inf\limits_{M}\psi_\varepsilon(t_2)&\geq&\sup\limits_{M}\psi_\varepsilon(t_2).
\end{eqnarray*}
By equation $(\ref{216})$, $\psi_\varepsilon(t)=\int_{0}^{t}e^{\beta(t-s)}f_\varepsilon(s)ds$. Hence $\tilde{\omega}_{\varepsilon}(t)=\omega_{\varepsilon}(t)$. \QEDB

\section{The long-time existence of the conical K\"ahler-Ricci flow with weak initial data}
\setcounter{equation}{0}

In this section, we study the long-time existence of the conical K\"ahler-Ricci flow $(\ref{TCKRF1})$ by the smooth approximation of the twisted K\"ahler-Ricci flows. We also prove the uniqueness of the conical K\"ahler-Ricci flow $(\ref{TCKRF1})$.

By Proposition \ref{200}, Lemma \ref{205} and Proposition \ref{210}, we conclude that for any $T>0$, there exists constants $C_1$ and $C_{2}$ depending only on $\|\varphi_0\|_{L^\infty(M)}$, $\beta$, $n$, $\omega_{0}$ and $T$, such that for any $\varepsilon>0$,
\begin{eqnarray}\label{219}&\ &\|\phi_{\varepsilon}(t)\|_{L^\infty([0,T]\times M)}\leq C_1,\\
\label{220}&\ &e^{-\frac{C_2}{t}}\omega_\varepsilon\leq\omega_{\varepsilon}(t)\leq e^{\frac{C_2}{t}}\omega_\varepsilon\ \ on\ (0,T]\times M.
\end{eqnarray}

We first prove the local uniform Calabi's $C^3$-estimate and curvature estimate along the flow $(\ref{TKRF2})$. Our proofs are similar as that in \cite{JWLXZ} (see section $2$ in \cite{JWLXZ} or section $3$ in \cite{MSBW}), but we need some arguments to handle the weak initial data case.

\begin{lem} For any $T>0$ and $B_r(p)\subset\subset M\setminus D$, there exist constants $C$, $C'$ and $C''$ such that for any $\varepsilon>0$ and $j\in\mathbb{N}^+$,
\begin{eqnarray*}S_{\varepsilon,j}&\leq&\frac{C'}{r^{2}}e^{\frac{C}{t}},\\
|Rm_{\varepsilon,j}|_{\omega_{\varepsilon,j}(t)}^{2}&\leq&\frac{C''}{r^{4}}e^{\frac{C}{t}}
\end{eqnarray*}
on $(0,T]\times B_{\frac{r}{2}}(p)$, where constants $C$, $C'$ and $C''$ depend only on $\|\varphi_0\|_{L^\infty(M)}$, $n$, $\beta$, $T$, $\omega_0$ and $dist_{\omega_0}(B_r(p),D)$.\end{lem}

{\bf Proof:}\ \ By Lemma \ref{205}, there exists uniform constat $C$ depending only on $\|\varphi_0\|_{L^\infty(M)}$, $n$, $\beta$, $T$, $\omega_0$ and $dist_{\omega_0}(B_r(p),D)$, such that
\begin{eqnarray}e^{-\frac{C}{t}}\omega_0\leq\omega_{\varepsilon,j}(t)\leq e^{\frac{C}{t}}\omega_0,\ \ on\ \ B_r(p)\times (0,T].
\end{eqnarray}

Let $r=r_0>r_1>\frac{r}{2}$ and $\psi$ be a nonnegative $C^\infty$ cut-off function that is identically equal to $1$ on $\overline{B_{r_{1}(p)}}$ and vanishes outside $B_r(p)$. We may assume that
\begin{eqnarray}|\partial\psi|_{\omega_{0}}^{2}\leq\frac{C}{r^{2}}\ \ \ and\ \ \ |\sqrt{-1}\partial\bar{\partial}\psi|_{\omega_{0}}\leq\frac{C}{r^{2}}.\end{eqnarray}
Straightforward calculations show that
\begin{eqnarray}(\frac{\partial}{\partial t}-\Delta_{\varepsilon,j})(\psi^2S_{\varepsilon,j})\leq \frac{C}{r^2} e^{\frac{C}{t}}S_{\varepsilon,j}+Ce^{\frac{C}{t}}.
\end{eqnarray}
By choosing sufficiently large $\alpha$, $\gamma$ and $A$, we get
 \begin{eqnarray*}&\ &(\frac{\partial}{\partial t}-\Delta_{\varepsilon,j})(e^{-\frac{2\alpha}{t}}\psi^2S_{\varepsilon,j}+Ae^{-\frac{2\gamma}{t}}tr_{\omega_0}\omega_{\varepsilon,j}(t))\\
 &\leq& -\frac{1}{r^2}e^{-\frac{3\gamma}{t}}S_{\varepsilon,j}+\frac{C}{r^2},
\end{eqnarray*}
where $\alpha=3\gamma$, $A=\frac{C+1}{r^{2}}$, constat $C$ depends only on $\|\varphi_0\|_{L^\infty(M)}$, $n$, $\beta$, $T$ $\omega_0$ and $dist_{\omega_{0}}(B_r(p),D)$. By the maximum principle, we conclude that
\begin{eqnarray*}S_{\varepsilon,j}&\leq&\frac{C'}{r^{2}}e^{\frac{6\gamma}{t}}\ \ \ on\ \ (0,T]\times B_{\frac{r}{2}}(p).
\end{eqnarray*}

Now we prove that $|Rm_{\varepsilon,j}|^{2}_{\omega_{\varepsilon,j}(t)}$ is uniformly bounded. Through computation, there exist uniform constants $C$ such that
\begin{eqnarray*}\label{3.20.9}&\ &(\frac{d}{dt}-\Delta_{\varepsilon,j})|Rm_{\varepsilon,j}|^{2}_{\omega_{\varepsilon,j}(t)}\\
&\leq& C|Rm_{\varepsilon,j}|^{3}_{\omega_{\varepsilon,j}(t)}+Ce^{\frac{C}{t}}|Rm_{\varepsilon,j}|^{2}_{\omega_{\varepsilon,j}(t)}+Ce^{\frac{C}{t}}|Rm_{\varepsilon,j}|_{\omega_{\varepsilon,j}(t)}+C e^{\frac{C}{t}}S_{\varepsilon,j}^{\frac{1}{2}}|Rm_{\varepsilon,j}|_{\omega_{\varepsilon,j}(t)}\\
&\ &+Ce^{\frac{C}{t}}S_{\varepsilon,j}|Rm_{\varepsilon,j}|_{\omega_{\varepsilon,j}(t)}-|\nabla_{\varepsilon,j}Rm_{\varepsilon,j}|^{2}_{\omega_{\varepsilon,j}(t)}-|\overline
{\nabla}_{\varepsilon,j}Rm_{\varepsilon,j}|_{\omega_{\varepsilon,j}(t)}^{2}+Ce^{\frac{C}{t}}\\ \nonumber
&\leq&C(|Rm_{\varepsilon,j}|^{3}_{\omega_{\varepsilon,j}(t)}+e^{\frac{\tau}{t}}+\frac{1}{r^{2}}e^{\frac{\tau}{t}}|Rm_{\varepsilon,j}|_{\omega_{\varepsilon,j}(t)})-|\nabla_{\varepsilon,j}
Rm_{\varepsilon,j}|_{\omega_{\varepsilon,j}(t)}^{2}-|\overline{\nabla}_{\varepsilon,j}Rm_{\varepsilon,j}|_{\omega_{\varepsilon,j}(t)}^{2}.
\end{eqnarray*}

Next, we show that $|Rm_{\varepsilon,j}|^{2}_{\omega_{\varepsilon,j}(t)}$ is uniformly bounded. We fix a smaller radius $r_{2}$ satisfying $r_{1}>r_{2}>\frac{r}{2}$. Let $\rho$ be a cut-off function identically equal to $1$ on $\overline{B_{r_{2}}}(p)$ and identically equal to $0$ outside $B_{r_{1}}$. We also let $\rho$ satisfy
$$|\partial\rho|_{\omega_{0}}^{2},\ |\sqrt{-1}\partial\bar{\partial}\rho|_{\omega_{0}}\leq\frac{C}{r^{2}}$$
for some uniform constant $C$. From the former part we know that $S_{\varepsilon,j}$ is bounded by $\frac{C}{r^{2}}e^{\frac{\tau}{t}}$ on $B_{r_{1}}(p)$. Let $K_t=\frac{\hat{C}}{r^{2}}e^{\frac{k\tau}{t}}$, $k$ and $\hat{C}$ be constants which are large enough such that $\frac{K_t}{2}\leq K_t-S_{\varepsilon,j}\leq K_t$. We consider
\begin{eqnarray}\label{3.20.11}F_{\varepsilon,j}=\rho^{2}e^{-\frac{2\delta}{t}}\frac{|Rm_{\varepsilon,j}|^{2}_{\omega_{\varepsilon,j}(t)}}{K_t-S_{\varepsilon,j}}
+Ae^{-\frac{2\sigma}{t}}S_{\varepsilon,j}.\end{eqnarray}
By computing, we have
\begin{eqnarray*}\label{3.20.12}&\ &(\frac{d}{dt}-\Delta_{\varepsilon,j})F_{\varepsilon,j}\\
&=&e^{-\frac{2\delta}{t}}\Big((-\triangle_{\varepsilon,j}\rho^{2})
\frac{|Rm_{\varepsilon,j}|^{2}_{\omega_{\varepsilon,j}(t)}}{K_t-S_{\varepsilon,j}}+\rho^{2}\frac{|Rm_{\varepsilon,j}|^{2}_{\omega_{\varepsilon,j}(t)}}
{(K_t-S_{\varepsilon,j})^{2}}(\frac{d}{dt}-\Delta_{\varepsilon,j})S_{\varepsilon,j}+\rho^{2}\frac{\hat{C}|Rm_{\varepsilon,j}|^{2}_{\omega_{\varepsilon,j}(t)}}
{r^2(K_t-S_{\varepsilon,j})^{2}}\frac{k\tau}{t^2}e^{\frac{k\tau}{t}}\\
&\ &+\rho^{2}\frac{1}{K_t-S_{\varepsilon,j}}(\frac{d}{dt}-\Delta_{\varepsilon,j})|Rm_{\varepsilon,j}|^{2}_{\omega_{\varepsilon,j}(t)}-4Re\langle\rho\frac{\nabla_{\varepsilon,j}\rho}{K_t-S_{\varepsilon,j}},\nabla_{\varepsilon,j}|Rm_{\varepsilon,j}|^{2}_{\omega_{\varepsilon,j}(t)}
\rangle_{\omega_{\varepsilon,j}(t)}\\
&\ &-4Re\langle\rho\frac{|Rm_{\varepsilon,j}|^{2}_{\omega_{\varepsilon,j}(t)}}{(K_t-S_{\varepsilon,j})^2}\nabla_{\varepsilon,j}S_{\varepsilon,j},\nabla_{\varepsilon,j}
\rho\rangle_{\omega_{\varepsilon,j}(t)}-2\frac{\rho^{2}|Rm_{\varepsilon,j}|^{2}_{\omega_{\varepsilon,j}(t)}}{(K_t-S_{\varepsilon,j})^3}|\nabla_{\varepsilon,j}S|
^{2}_{\omega_{\varepsilon,j}(t)}\\
&\ &-2Re\langle\rho^{2}\frac{\nabla_{\varepsilon,j}S_{\varepsilon,j}}{(K_t-S_{\varepsilon,j})^2},\nabla_{\varepsilon,j}|
Rm_{\varepsilon,j}|^{2}_{\omega_{\varepsilon,j}(t)}
\rangle_{\omega_{\varepsilon,j}(t)}\Big)+Ae^{-\frac{2\sigma}{t}}(\frac{d}{dt}-\triangle_{\omega_{\varepsilon,j}(t)})S_{\varepsilon,j}\\
&\ &+A\frac{2\sigma}{t^2}e^{-\frac{2\sigma}{t}}S_{\varepsilon,j}+\frac{2\delta}{t^2}e^{-\frac{2\delta}{t}}\rho^{2}\frac{|Rm_{\varepsilon,j}|^{2}_
{\omega_{\varepsilon,j}(t)}}{K_t-S_{\varepsilon,j}}.
\end{eqnarray*}

We only consider an inner point $(t_{0},x_{0})$ which is a maximum point of $F_{\varepsilon,j}$ achieved on $[0,T]\times\overline{B_{r_{1}}(p)}$. We use the fact that $\nabla_{\varepsilon,j} F_{\varepsilon,j}=0$ at this point, 
\begin{eqnarray*}\label{3.20.13} &\ &e^{-\frac{2\delta}{t}}\Big(2\rho\nabla_{\varepsilon,j}\rho\frac{|Rm_{\varepsilon,j}|^{2}_{\omega_{\varepsilon,j}(t)}}{K_t-S_{\varepsilon,j}}+\rho^2\frac{\nabla_{\varepsilon,j}
|Rm_{\varepsilon,j}|^{2}
_{\omega_{\varepsilon,j}(t)}}{K_t-S_{\varepsilon,j}}+\rho^2\frac{|Rm_{\varepsilon,j}|^{2}_{\omega_{\varepsilon,j}(t)}
\nabla_{\varepsilon,j}S}{(K_t-S_{\varepsilon,j})^2}\Big)\\
&\ &+Ae^{-\frac{2\sigma}{t}}\nabla_{\varepsilon,
j}S_{\varepsilon,j}=0.\end{eqnarray*}
Combining the above two equalities, we have

\begin{eqnarray*}\label{3.20.12}&\ &(\frac{d}{dt}-\Delta_{\varepsilon,j})F_{\varepsilon,j}\\
&=&e^{-\frac{2\delta}{t}}\Big((-\triangle_{\varepsilon,j}\rho^{2})
\frac{|Rm_{\varepsilon,j}|^{2}_{\omega_{\varepsilon,j}(t)}}{K_t-S_{\varepsilon,j}}+\rho^{2}\frac{|Rm_{\varepsilon,j}|^{2}_{\omega_{\varepsilon,j}(t)}}
{(K_t-S_{\varepsilon,j})^{2}}(\frac{d}{dt}-\Delta_{\varepsilon,j})S_{\varepsilon,j}+\rho^{2}\frac{\hat{C}|Rm_{\varepsilon,j}|^{2}_{\omega_{\varepsilon,j}(t)}}
{r^2(K_t-S_{\varepsilon,j})^{2}}\frac{k\tau}{t^2}e^{\frac{k\tau}{t}}\\
&\ &+\rho^{2}\frac{1}{K_t-S_{\varepsilon,j}}(\frac{d}{dt}-\Delta_{\varepsilon,j})|Rm_{\varepsilon,j}|^{2}_{\omega_{\varepsilon,j}(t)}
-4Re\langle\rho\frac{\nabla_{\varepsilon,j}\rho}{K_t-S_{\varepsilon,j}},\nabla_{\varepsilon,j}|Rm_{\varepsilon,j}|^{2}_{\omega_{\varepsilon,j}(t)}
\rangle_{\omega_{\varepsilon,j}(t)}\Big) \\
&\ & +2Ae^{-\frac{2\sigma}{t}}\frac{|\nabla_{\varepsilon,j}S_{\varepsilon,j}|^{2}_{\omega_{\varepsilon,j}(t)}}{K_t-S_{\varepsilon,j}}+Ae^{-\frac{2\sigma}{t}}
(\frac{d}{dt}-\triangle_{\omega_{\varepsilon,j}(t)})S_{\varepsilon,j}+A\frac{2\sigma}{t^2}e^{-\frac{2\sigma}{t}}S_{\varepsilon,j}\\
&\ &+\frac{2\delta}{t^2}e^{-\frac{2\delta}{t}}\rho^{2}\frac{|Rm_{\varepsilon,j}|^{2}_
{\omega_{\varepsilon,j}(t)}}{K_t-S_{\varepsilon,j}}.
\end{eqnarray*}

Our goal is to show that at $(t_{0},x_{0})$ we have $e^{-\frac{2\delta}{t}}|Rm_{\varepsilon,j}|^{2}_{\omega_{\varepsilon,j}(t)}\leq\frac{C}{r^{4}}$ for some uniform constant $C$ and $\delta$. Without loss of generality, we assume that $|Rm_{\varepsilon,j}|^{3}_{\omega_{\varepsilon,j}(t)}\geq e^{\frac{\tau}{t}}+\frac{1}{r^{2}}e^{\frac{\tau}{t}}|Rm_{\varepsilon,j}|_{\omega_{\varepsilon,j}(t)}$ at $(t_{0},x_{0})$.
\begin{eqnarray*}\label{3.20.15}(\frac{d}{dt}-\triangle_{\omega_{\varepsilon,j}(t)})|Rm_{\varepsilon,j}|^{2}_{\omega_{\varepsilon,j}(t)}&\leq & C|Rm_{\varepsilon,j}|^{3}_{\omega_{\varepsilon,j}(t)}-|\nabla_{\varepsilon,j}Rm_{\varepsilon,j}|_{\omega_{\varepsilon,j}(t)}^{2}-|\overline{\nabla}_{\varphi}Rm_{\varepsilon,j}|_{\omega_{\varepsilon,j}(t)}^{2},\\
\label{3.20.17}(\frac{d}{dt}-\triangle_{\omega_{\varepsilon,j}(t)})S_{\varepsilon,j}&\leq & \frac{C}{r^{2}}e^{\frac{\tau}{t}}-|\nabla_{\varepsilon,j}X|_{\omega_{\varepsilon,j}(t)}^{2}-|\overline{\nabla}_{\varepsilon,j}X|_{\omega_{\varepsilon,j}(t)}^{2}
\end{eqnarray*}
on $B_{r_{1}}(p)$.
We also note that
\begin{eqnarray*}\label{3.20.16}|\nabla_{\varepsilon,j}|Rm_{\varepsilon,j}|_{\omega_{\varepsilon,j}(t)}^{2}|_{\omega_{\varepsilon,j}(t)}&\leq&|Rm_{\varepsilon,j}|_{\omega_{\varepsilon,j}(t)}(|\nabla_{\varepsilon,j}Rm_{\varepsilon,j}|_{\omega_{\varepsilon,j}(t)}+|\overline{\nabla}_{\varepsilon,j}Rm_{\varepsilon,j}|_{\omega_{\varepsilon,j}(t)}),\\
\label{3.20.16'}|\nabla_{\varepsilon,j}S_{\varepsilon,j}|^{2}_{\omega_{\varepsilon,j}(t)}&\leq&2S_{\varepsilon,j}(|\nabla_{\varepsilon,j}X|_{\omega_{\varepsilon,j}(t)}^{2}+|\overline{\nabla}_{\varepsilon,j}X|_{\omega_{\varepsilon,j}(t)}^{2}).
\end{eqnarray*}
By using the above inequalities, at $(t_{0},x_{0})$, we have
\begin{eqnarray*}\label{3.20.18}&\ &(\frac{d}{dt}-\triangle_{\omega_{\varepsilon,j}(t)})F_{\varepsilon,j}\\
&\leq&-Ae^{-\frac{2\sigma}{t}}(|\nabla_{\varepsilon,j}X|_{\omega_{\varepsilon,j}(t)}^{2}+|\overline{\nabla}_{\varepsilon,j}X|_{\omega_{\varepsilon,j}(t)}^{2})
+\frac{AC}{r^{2}}e^{-\frac{2\sigma}{t}}e^{\frac{\tau}{t}}+e^{-\frac{2\delta}{t}}\Big(\frac{Ce^{\frac{C}{t}}|Rm_{\varepsilon,j}|_{\omega_{\varepsilon,j}(t)}^{2}}{K_tr^{2}}
\\
&\ &-\frac{\rho^{2}|Rm_{\varepsilon,j}|_{\omega_{\varepsilon,j}(t)}^{2}(|\nabla_{\varepsilon,j}X|_{\omega_{\varepsilon,j}(t)}^{2}
+|\overline{\nabla}_{\varepsilon,j}X|_{\omega_{\varepsilon,j}(t)}^{2})}{K_t^{2}}
+\frac{C\rho^{2}|Rm_{\varepsilon,j}|_{\omega_{\varepsilon,j}(t)}^{3}}{K_t}\\
&\ &-\frac{\rho^{2}(|\nabla_{\varepsilon,j}Rm_{\varepsilon,j}|_{\omega_{\varepsilon,j}(t)}^{2}+|\overline{\nabla}_{\varepsilon,j}Rm_{\varepsilon,j}|
_{\omega_{\varepsilon,j}(t)}^{2})}{K_t}+\frac{Ce^{\frac{C}{t}}|Rm_{\varepsilon,j}|_{\omega_{\varepsilon,j}(t)}^{2}}{K_tr^{2}}+\frac{C\rho^{2}|Rm_{\varepsilon,j}|
_{\omega_{\varepsilon,j}(t)}^{2}}{K_t^{2}r^{2}}e^{\frac{\tau}{t}}\\
&\ &+\frac{\rho^{2}(|\nabla_{\varepsilon,j}Rm_{\varepsilon,j}|_{\omega_{\varepsilon,j}(t)}^{2}+|\overline{\nabla}_{\varepsilon,j}Rm_{\varepsilon,j}|
_{\omega_{\varepsilon,j}(t)}^{2})}{K_t}\Big)+\frac{8Ae^{-\frac{2\sigma}{t}}S_{\varepsilon,j}(|\nabla_{\varepsilon,j}X|_{\omega_{\varepsilon,j}(t)}^{2}
+|\overline{\nabla}_{\varepsilon,j}X|_{\omega_{\varepsilon,j}(t)}^{2})}{K_t}\\
&\ &+\rho^{2}e^{-\frac{2\delta}{t}}\frac{2\hat{C}|Rm_{\varepsilon,j}|^{2}_{\omega_{\varepsilon,j}(t)}}
{r^2K_t^{2}}\frac{k\tau}{t^2}e^{\frac{k\tau}{t}}+A\frac{2\sigma}{t^2}e^{-\frac{2\sigma}{t}}S_{\varepsilon,j}+\frac{4\delta}{t^2}e^{-\frac{2\delta}{t}}\rho^{2}\frac{|Rm_{\varepsilon,j}|^{2}_
{\omega_{\varepsilon,j}(t)}}{K_t}.
\end{eqnarray*}
Let $\hat{C}$ be sufficiently large so that $\frac{8AS_{\varepsilon,j}Q}{K_t}\leq
\frac{AQ}{2}$, where we denote $Q=|\nabla_{\varepsilon,j}X|_{\omega_{\varepsilon,j}(t)}^{2}+|\overline{\nabla}_{\varepsilon,j}X|_{\omega_{\varepsilon,j}(t)}^{2}$. Then

\begin{eqnarray}\label{3.20.19}\nonumber\frac{C\rho^{2}|Rm_{\varepsilon,j}|_{\omega_{\varepsilon,j}(t)}^{3}}{K_t}&\leq&\frac{\rho^{2}|Rm_{\varepsilon,j}|
_{\omega_{\varepsilon,j}(t)}^{4}}{2K_t^{2}}+C\rho^{2}|Rm_{\varepsilon,j}|_{\omega_{\varepsilon,j}(t)}^{2}\\
&\leq&\frac{\rho^{2}|Rm_{\varepsilon,j}|_{\omega_{\varepsilon,j}(t)}^{2}Q}{K_t^{2}}+Ce^{\frac{C}{t}}\rho^{2}|Rm_{\varepsilon,j}|_{\omega_{\varepsilon,j}(t)}^{2}
\end{eqnarray}
Let $k=1$, $\delta=2\sigma$ and $\tau-2\sigma<0$, where $\sigma$ is sufficiently large. We conclude that the evolution equation of $F_{\varepsilon,j}$ can be controlled as follows,
\begin{eqnarray*}\label{3.20.20}(\frac{d}{dt}-\Delta_{\varepsilon,j})F_{\varepsilon,j}&\leq&-\frac{Ae^{-\frac{2\sigma}{t}}Q}{2}+\frac{AC}{r^{2}}
e^{-\frac{2\sigma}{t}}e^{\frac{\tau}{t}}+\frac{AC}{t^2r^{2}}
e^{-\frac{2\sigma}{t}}e^{\frac{\tau}{t}}+Ce^{-\frac{\delta}{t}}|Rm_{\varepsilon,j}|_{\omega_{\varepsilon,j}(t)}^{2}\\ \nonumber
&\leq&-\frac{Ae^{-\frac{2\sigma}{t}}Q}{2}+\frac{AC}{r^{2}}e^{-\frac{2\sigma}{t}}e^{\frac{\tau}{t}}+\frac{AC}{t^2r^{2}}
e^{-\frac{2\sigma}{t}}e^{\frac{\tau}{t}}+Ce^{-\frac{\delta}{t}}Q+Ce^{-\frac{\delta}{2t}}.
\end{eqnarray*}
Now we choose a sufficiently large $A$ such that $A=2(C+1)$ and obtain
$$e^{-\frac{\delta}{t}}Q\leq\frac{C}{r^{2}}$$
at $(t_{0},x_{0})$. This implies that $e^{-\frac{2\delta}{t}}|Rm_{\varepsilon,j}|_{\omega_{\varepsilon,j}(t)}^{2}\leq\frac{C}{r^{2}}$ at this point, where $C$ depends only on $\|\varphi_0\|_{L^\infty(M)}$, $n$, $\beta$, $T$, $dist_{\omega_{0}}(B_r(p),D)$, $\|\theta\|_{C^2(B_{r}(p))}$ and $\omega_0$. Following that we conclude that $F_{\varepsilon,j}$ is bounded by $\frac{C}{r^{2}}$ at $(t_{0},x_{0})$. Hence on $[0,T]\times \overline{B_{r_{2}}(p)}$, we obtain
\begin{eqnarray}\label{3.20.21}|Rm_{\varepsilon,j}|_{\omega_{\varepsilon,j}(t)}^{2}\leq\frac{C}{r^{4}}e^{\frac{2\delta+\tau}{t}},\end{eqnarray}
where $C$, $\delta$ and $\tau$ depend only on $\|\varphi_0\|_{L^\infty(M)}$, $n$, $\beta$, $T$, $dist_{\omega_{0}}(B_r(p),D)$ and $\omega_0$.\QEDB

\medskip

By using the standard parabolic Schauder regularity theory \cite{GLIE}, we obtain the following proposition.

\begin{pro}\label{217} For any $0<\delta<T<\infty$, $k\in\mathbb{N}^{+}$ and $B_r(p)\subset\subset M\setminus D$, there exists constant $C_{\delta,T,k,p,r}$ depends only on $\|\varphi_0\|_{L^\infty(M)}$, $n$, $\beta$, $\delta$, $k$, $T$, $dist_{\omega_{0}}(B_r(p),D)$ and $\omega_0$, such that for any $\varepsilon>0$ and $j\in\mathbb{N}^+$,
\begin{eqnarray}\|\varphi_{\varepsilon,j}(t)\|_{C^{k}\big([\delta,T]\times B_r(p)\big)}\leq C_{\delta,T,k,p,r}.
\end{eqnarray}
\end{pro}

Through a further observation to equation $(\ref{CMAE3})$, we prove the monotonicity of $\varphi_\varepsilon(t)$ with respect to $\varepsilon$.

\begin{pro}\label{218} For any $(t,x)\in [0,T]\times M$, $\varphi_\varepsilon(t,x)$ is monotone decreasing as $\varepsilon\searrow0$.
\end{pro}

{\bf Proof:}\ \ For any $\varepsilon_1<\varepsilon_2$, let $\psi_{1,2}(t)=\varphi_{\varepsilon_1}(t)-\varphi_{\varepsilon_2}(t)$. Then we have

\begin{eqnarray}\nonumber&\ &\frac{\partial }{\partial t}(e^{-\beta t}\psi_{1,2}(t))\\
&\leq&e^{-\beta t}\log\frac{\big(e^{-\beta t}\omega_{0}+\sqrt{-1}\partial\bar{\partial}e^{-\beta t}\varphi_{\varepsilon_2}(t)+\sqrt{-1}\partial\bar{\partial}e^{-\beta t}\psi_{1,2}(t)\big)^{n}}
  {(e^{-\beta t}\omega_{0}+\sqrt{-1}\partial\bar{\partial}e^{-\beta t}\varphi_{\varepsilon_2}(t))^{n}}.
\end{eqnarray}

Let $\tilde{\psi}_{1,2}(t)=e^{-\beta t}\psi_{1,2}(t)-\delta t$ with $\delta>0$ and $(t_0, x_0)$ be the maximum point of $\tilde{\psi}_{1,2}(t)$ on $[0,T]\times M$. If $t_0>0$, by maximum principle, at this point, we have
\begin{eqnarray}0\leq\frac{\partial }{\partial t}\tilde{\psi}_{1,2}(t)=\frac{\partial }{\partial t}(e^{-\beta t}\psi_{1,2}(t))-\delta\leq-\delta
\end{eqnarray}
which is impossible, hence $t_0=0$. So for any $(t,x)\in [0,T]\times M$,
\begin{eqnarray}\psi_{1,2}(t,x)\leq e^{\beta t}\sup\limits_{M}\psi_{1,2}(0,x)+Te^{\beta T}\delta=Te^{\beta T}\delta.
\end{eqnarray}
Let $\delta\rightarrow0$, we conclude that $\varphi_{\varepsilon_1}(t,x)\leq\varphi_{\varepsilon_2}(t,x)$.\QEDB

\medskip

For any $[\delta, T]\times K\subset\subset (0,\infty)\times M\setminus D$ and $k\geq0$, $\|\varphi_{\varepsilon,j}(t)\|_{C^{k}([\delta,T]\times K)}$ is uniformly bounded by Proposition \ref{217}. Let $j$ approximate to $\infty$, we obtain that $\|\varphi_{\varepsilon}(t)\|_{C^{k}([\delta,T]\times K)}$ is uniformly bounded. Then let $\delta$ approximate to $0$, $T$ approximate to $\infty$ and $K$ approximate to $M\setminus D$, by diagonal rule, we get a sequence $\{\varepsilon_i\}$, such that $\varphi_{\varepsilon_i}(t)$ converges in $C^\infty_{loc}$ topology on $(0,\infty)\times (M\setminus D)$ to a function $\varphi(t)$ that is smooth on $C^\infty\big((0,\infty)\times(M\setminus D)\big)$ and satisfies equation
\begin{eqnarray}\label{223}  \frac{\partial \varphi(t)}{\partial t}=\log\frac{(\omega_{0}+\sqrt{-1}\partial\bar{\partial}\varphi(t))^{n}}{\omega_{0}^{n}}+F_{0}+\beta\varphi(t)
  +\log|s|_{h}^{2(1-\beta)}
\end{eqnarray}
on $(0,\infty)\times (M\setminus D)$. Since $\varphi_\varepsilon(t)$ is monotone decreasing as $\varepsilon\rightarrow0$, we conclude that $\varphi_{\varepsilon}(t)$ converges in $C^\infty_{loc}$ topology on $(0,\infty)\times (M\setminus D)$ to $\varphi(t)$. Combining the above arguments with
$(\ref{219})$ and $(\ref{220})$, for any $T>0$, we have
\begin{eqnarray}\label{221}&\ &\|\varphi(t)\|_{L^\infty\big((0,T]\times (M\setminus D)\big)}\leq C_1,\\
\label{222}&\ &e^{-\frac{C_2}{t}}\omega_\beta\leq\omega(t)\leq e^{\frac{C_2}{t}}\omega_\beta\ \ on\ (0,T]\times (M\setminus D),
\end{eqnarray}
where $\omega(t)=\omega_0+\sqrt{-1}\partial\overline{\partial}\varphi(t)$, constants $C_1$ and $C_{2}$ depend only on $\|\varphi_0\|_{L^\infty(M)}$, $\beta$, $n$, $\omega_{0}$ and $T$.

\medskip

\begin{pro}\label{2188} For any $t>0$, $\varphi(t)$ is H\"older continuous on $M$ with respect to the metric $\omega_0$.
\end{pro}

\medskip

{\bf Proof:}\ \ We assume that $t\in[\delta, T]$ for some $\delta$ and $T$ satisfying $0<\delta<T<\infty$. By $(\ref{222})$ we have
\begin{eqnarray}C^{-1}\omega_\beta\leq\omega(t)\leq C\omega_\beta\ \ on\ [\delta,T]\times  (M\setminus D),
\end{eqnarray}
where constant $C$ depends only on $\|\varphi_0\|_{L^\infty(M)}$, $\beta$, $T$, $n$, $\omega_{0}$ and $\delta$. Combining this estimate and the fact that
$\log\frac{\omega_\beta^n|s|_h^{2(1-\beta)}}{\omega_0^n}$ is bounded uniformly on $M\setminus D$, we obtain
\begin{eqnarray}\|\log\frac{\omega^{n}(t)|s|_h^{2(1-\beta)}}{\omega_0^n}\|_{L^\infty\big([\delta,T]\times (M\setminus D)\big)}\leq C
\end{eqnarray}
for some uniform constant $C$ independent of $t$. Therefore, $\|\frac{\partial \varphi(t)}{\partial t}\|_{L^\infty\big([\delta,T]\times (M\setminus D)\big)}$ is uniformly bounded by equation $(\ref{223})$ and estimate $(\ref{221})$. We rewrite equation $(\ref{223})$ as
\begin{eqnarray}\label{224}(\omega_{0}+\sqrt{-1}\partial\bar{\partial}\varphi(t))^{n}=e^{\frac{\partial \varphi(t)}{\partial t}-F_{0}-\beta\varphi(t)}\frac{\omega_{0}^{n}}{|s|_{h}^{2(1-\beta)}}.
\end{eqnarray}
The function on the right side of equation $(\ref{224})$ is $L^p$ integrable with respect to $\omega_0^{n}$ for some $p>1$. By S. Kolodziej's $L^p$-estimates \cite{K}, we know that $\varphi(t)$ is H\"older continuous on $M$ with respect to $\omega_0$ for any $t>0$.\QEDB

\medskip

Next, by using the monotonicity of $\varphi_{\varepsilon}(t)$ with respect to $\varepsilon$ and constructing auxiliary function, we prove the continuity of $\varphi(t)$ as $t\rightarrow0^{+}$.

\begin{pro}\label{101} $\varphi(t)\in C^0([0,\infty)\times M)$ and
\begin{eqnarray}\label{102}\lim\limits_{t\rightarrow0^+}\|\varphi(t)-\varphi_{0}\|_{L^{\infty}(M)}=0.
\end{eqnarray}
\end{pro}

{\bf Proof:}\ \ Through the above arguments, we only need prove limit $(\ref{102})$. By the monotonicity of $\varphi_{\varepsilon}(t)$ with respect to $\varepsilon$, for any $(t,z)\in(0,T]\times M$, we have
\begin{eqnarray}\label{2210}\nonumber\varphi(t,z)-\varphi_0(z)&\leq&\varphi_{\varepsilon_{1}}(t,z)-\varphi_0(z)\\\nonumber
&\leq&|\varphi_{\varepsilon_1}(t,z)-\varphi_{\varepsilon_1,j}(t,z)|
+|\varphi_{\varepsilon_1,j}(t,z)-\varphi_{0,j}(z)|\\
&\ &+|\varphi_{0,j}(z)-\varphi_{0}(z)|.
\end{eqnarray}
Since $\varphi_{\varepsilon_1,j}(t)$ is a Cauchy sequence in $L^\infty([0,T]\times M)$,
\begin{eqnarray}\lim\limits_{j\rightarrow\infty}\|\varphi_{\varepsilon_1}(t,z)-\varphi_{\varepsilon_1,j}(t,z)\|_{L^{\infty}([0,T]\times M)}=0.
\end{eqnarray}
From $(\ref{00000})$, we have
\begin{eqnarray}
\lim\limits_{j\rightarrow\infty}\|\varphi_{0,j}(z)-\varphi_{0}(z)\|_{L^{\infty}(M)}=0,
\end{eqnarray}
For any $\epsilon>0$, there exists $N$ such that for any $j>N$,
\begin{eqnarray*}\sup\limits_{[0,T]\times M}|\varphi_{\varepsilon_1}(t,z)-\varphi_{\varepsilon_1,j}(t,z)|&<&\frac{\epsilon}{3},\\
\sup\limits_{M}|\varphi_{0,j}(z)-\varphi_{0}(z)|&<&\frac{\epsilon}{3}.
\end{eqnarray*}
Fix such $\varepsilon_{1}$ and $j$, there exists $0<\delta_1<T$ such that
\begin{eqnarray}\sup\limits_{[0,\delta_1]\times M}|\varphi_{\varepsilon_1,j}(t,z)-\varphi_{0,j}|<\frac{\epsilon}{3}.
\end{eqnarray}
Combining the above estimates together, for any $t\in(0,\delta_1]$ and $z\in M$,
\begin{eqnarray}\varphi(t,z)-\varphi_{0}(z)<\epsilon.
\end{eqnarray}

On the other hand, by S. Kolodziej's results \cite{K000}, there exists a smooth solution $u_{\varepsilon,j}$ to the equation
\begin{eqnarray}(\omega_{0}+\sqrt{-1}\partial\bar{\partial}u_{\varepsilon,j})^{n}=e^{-F_{0}-\beta\varphi_{0,j}+\hat{C}}\frac{\omega_{0}^{n}}
{(\varepsilon^2+|s|_{h}^{2})^{(1-\beta)}},
\end{eqnarray}
and $u_{\varepsilon,j}$ satisfies \begin{eqnarray}\|u_{\varepsilon,j}\|_{L^\infty(M)}\leq C,
\end{eqnarray}
where $\hat{C}$ is a uniform normalization constant independent of $\varepsilon$ and $j$, constant $C$ depends only on $\|\varphi_0\|_{L^\infty(M)}$, $\beta$ and $F_0$.

We define function
\begin{eqnarray}\psi_{\varepsilon,j}(t)=(1-te^{\beta t})\varphi_{0,j}+te^{\beta t}u_{\varepsilon,j}+h(t)e^{\beta t},
\end{eqnarray}
where
\begin{eqnarray}\nonumber h(t)&=&-t\|\varphi_{0,j}\|_{L^\infty(M)}-t\|u_{\varepsilon,j}\|_{L^\infty(M)}+n(t\log t-t)e^{-\beta t}\\
&\ &+\beta n\int_0^te^{-\beta s}s\log s ds-\frac{\hat{C}}{\beta}e^{-\beta t}+\frac{\hat{C}}{\beta}
\end{eqnarray}
and $h(0)=0$. 

Straightforward calculations show that
\begin{eqnarray*}\frac{\partial}{\partial t}\psi_{\varepsilon,j}(t)-\beta\psi_{\varepsilon,j}(t)&=&-\beta \varphi_{0,j}-e^{\beta t}\varphi_{0,j}+e^{\beta t}u_{\varepsilon,j}+e^{\beta t}\frac{\partial}{\partial t}h(t)\\
&=&-\beta \varphi_{0,j}-e^{\beta t}\varphi_{0,j}+e^{\beta t}u_{\varepsilon,j}-e^{\beta t}\|\varphi_{0,j}\|_{L^\infty(M)}-e^{\beta t}\|u_{\varepsilon,j}\|_{L^\infty(M)}\\
&\ &+n\log t-\beta n(t\log t-t)+\beta nt\log t+\hat{C}\\
&\leq& -\beta \varphi_{0,j}+n\log t+n\beta t+\hat{C}.
\end{eqnarray*}
Therefore, we have
\begin{eqnarray*}e^{\frac{\partial}{\partial t}\psi_{\varepsilon,j}(t)-\beta\psi_{\varepsilon,j}(t)}\omega_0^n\leq
t^ne^{n\beta t}e^{-\beta \varphi_{0,j}+\hat{C}}\omega_0^n.
\end{eqnarray*}
When $t$ is sufficiently small,
\begin{eqnarray*}\omega_0+\sqrt{-1}\partial\overline{\partial}\psi_{\varepsilon,j}(t)&=&(1-te^{\beta t})(\omega_0+\sqrt{-1}\partial\overline{\partial}\varphi_{0,j})+te^{\beta t}(\omega_0+\sqrt{-1}\partial\overline{\partial}u_{\varepsilon,j})\\
&\geq&te^{\beta t}(\omega_0+\sqrt{-1}\partial\overline{\partial}u_{\varepsilon,j}).
\end{eqnarray*}
Combining the above inequalities,
\begin{eqnarray*}(\omega_0+\sqrt{-1}\partial\overline{\partial}\psi_{\varepsilon,j}(t))^n&\geq&t^ne^{n\beta t}(\omega_0+\sqrt{-1}\partial\overline{\partial}u_{\varepsilon,j})^n\\
&=&t^ne^{n\beta t}e^{-F_{0}-\beta\varphi_{0,j}+\hat{C}}\frac{\omega_{0}^{n}}
{(\varepsilon^2+|s|_{h}^{2})^{(1-\beta)}}\\
&\geq&e^{-F_{0}+\frac{\partial}{\partial t}\psi_{\varepsilon,j}(t)-\beta\psi_{\varepsilon,j}(t)}\frac{\omega_{0}^{n}}
{(\varepsilon^2+|s|_{h}^{2})^{(1-\beta)}}.
\end{eqnarray*}
This equation is equivalent to
\begin{eqnarray}
\begin{cases}
  \frac{\partial}{\partial t}\psi_{\varepsilon,j}(t)\leq \log \frac{(\omega_0+\sqrt{-1}\partial\overline{\partial}\psi_{\varepsilon,j}(t))^n}{\omega_{0}^{n}}+\beta\psi_{\varepsilon,j}(t)\\
  \ \ \ \ \ \ \ \ \ \ \ \ \ \ +F_{0}+\log
(\varepsilon^2+|s|_{h}^{2})^{(1-\beta)}.\\
  \psi_{\varepsilon,j}(0)=\varphi_{0,j}\\
  \end{cases}
\end{eqnarray}

Let $\tilde{\psi}_{\varepsilon,j}(t)=\varphi_{\varepsilon,j}(t)-\psi_{\varepsilon,j}(t)$, then
\begin{eqnarray}
\begin{cases}
  \frac{\partial}{\partial t}\tilde{\psi}_{\varepsilon,j}(t)\geq \log
   \frac{(\omega_0+\sqrt{-1}\partial\overline{\partial}\psi_{\varepsilon,j}(t)+\sqrt{-1}\partial\overline{\partial}\tilde{\psi}_{\varepsilon,j}(t))^n}
  {(\omega_0+\sqrt{-1}\partial\overline{\partial}\psi_{\varepsilon,j}(t))^n}+\beta\tilde{\psi}_{\varepsilon,j}(t).\\
  \\
  \tilde{\psi}_{\varepsilon,j}(0)=0\\
  \end{cases}
\end{eqnarray}
By the similar arguments as that in the proof of Proposition \ref{218}, for any $(t,z)\in[0,T]\times M$,
\begin{eqnarray}\tilde{\psi}_{\varepsilon,j}(t,z)\geq0,
\end{eqnarray}
That is,  for any $(t,z)\in[0,T]\times M$
\begin{eqnarray}\nonumber\varphi_{\varepsilon,j}(t,z)-\varphi_{0,j}(z)&\geq&-te^{\beta t}\varphi_{0,j}+te^{\beta t}u_{\varepsilon,j}+h(t)e^{\beta t}\\
&\geq&-C(e^{\beta t}-1)+h_{1}(t)e^{\beta t} ,
\end{eqnarray}
where $h_{1}(t)=-n(t\log t-t)e^{-\beta t}+\beta n\int_0^te^{-\beta s}s\log s ds$, constant $C$ depends only on $\|\varphi_0\|_{L^\infty(M)}$, $\beta$ and $F_0$. Let $j\rightarrow\infty$ and then $\varepsilon\rightarrow0$, we have
\begin{eqnarray}\varphi(t,z)-\varphi_0(z)\geq-C(e^{\beta t}-1)+h_{1}(t)e^{\beta t}.
\end{eqnarray}
There exists $\delta_2$ such that for any $t\in[0,\delta_2]$,
\begin{eqnarray}-C(e^{\beta t}-1)+h_{1}(t)e^{\beta t}>-\epsilon.
\end{eqnarray}
Let $\delta=\min(\delta_1,\delta_2)$, then for any $t\in(0,\delta]$ and $z\in M$,
\begin{eqnarray}-\epsilon<\varphi(t,z)-\varphi_{0}(z)<\epsilon.
\end{eqnarray}
This completes the proof of the proposition.\QEDB

\begin{thm}\label{225} $\omega(t)=\omega_{0}+\sqrt{-1}\partial\bar{\partial}\varphi(t)$ is a long-time solution to the conical K\"ahler-Ricci flow $(\ref{TCKRF1})$.
\end{thm}

{\bf Proof:}\ \ We should only prove that $\omega(t)$ satisfies equation $(\ref{TCKRF1})$ in the sense of currents on $[0,\infty)\times M$.

Let $\eta=\eta(t,x)$ be a smooth $(n-1,n-1)$-form with compact support in $(0,\infty)\times M$. Without loss of generality, we assume that its compact support included in $(\delta,T)$ ($0<\delta<T<\infty$). On $[\delta,T]\times M$, by $(\ref{219})$, $(\ref{220})$, $(\ref{221})$ and $(\ref{222})$, we know that $\log\frac{\omega_{\varepsilon}^{n}(t)(\varepsilon^{2}+|s|_{h}^{2})^{1-\beta}}{\omega_{0}^{n}}$, $\log\frac{\omega^{n}(t)|s|_{h}^{2(1-\beta)}}{\omega_{0}^{n}}$, $\varphi_{\varepsilon}(t)$ and $\varphi(t)$ are uniformly bounded by constants depending only on $\|\varphi_0\|_{L^\infty(M)}$, $n$, $\beta$, $\delta$ and $T$. On $[\delta,T]$, we have
\begin{eqnarray}\label{201503201}\nonumber&\ &\int_{M}\frac{\partial\omega_{\varepsilon}(t)}{\partial t}\wedge\eta=\int_{M}\sqrt{-1}\partial\bar{\partial}\frac{\partial\varphi_{\varepsilon}(t)}{\partial t}\wedge\eta\\\nonumber
&=&\int_{M}\sqrt{-1}\partial\bar{\partial}\big(\log\frac{\omega_{\varepsilon}^{n}(\varepsilon^{2}+|s|_{h}^{2})^{1-\beta}}
{\omega_{0}^{n}}+F_{0}+\beta\varphi_{\varepsilon}(t)\big)\wedge\eta\\\nonumber
&=&\int_{M}\log\big(\log\frac{\omega_{\varepsilon}^{n}(\varepsilon^{2}+|s|_{h}^{2})^{1-\beta}}
{\omega_{0}^{n}}+F_{0}+\beta\varphi_{\varepsilon}(t)\big)
  \sqrt{-1}\partial\bar{\partial}\eta\\\nonumber
&\xrightarrow{\varepsilon\rightarrow 0}&\int_{M}(\log\frac{\omega^{n}(t)}{\omega_{0}^{n}}+F_{0}+\beta\varphi(t)
  +\log|s|_{h}^{2(1-\beta)})\sqrt{-1}\partial\bar{\partial}\eta\\\nonumber
&=&\int_{M}\sqrt{-1}\partial\bar{\partial}(\log\frac{\omega_{\varphi(t)}^{n}}{\omega_{0}^{n}}+F_{0}+\beta\varphi(t)
  +\log|s|_{h}^{2(1-\beta)})\wedge\eta\\
&=&\int_{M} (-Ric (\omega(t))+\beta \omega(t) +2\pi (1-\beta ) [D])\wedge\eta.
\end{eqnarray}
At the same time, there also holds
\begin{eqnarray}\label{201503202}\nonumber\int_M \omega_{\varphi_{\varepsilon}(t)}\wedge\frac{\partial\eta}{\partial t}&=&\int_M \omega_{0}\wedge\frac{\partial\eta}{\partial t}+\int_M \sqrt{-1}\partial\bar{\partial}\varphi_{\varepsilon}(t)\wedge\frac{\partial\eta}{\partial t}\\\nonumber
&=&\int_M \omega_{0}\wedge\frac{\partial\eta}{\partial t}+\int_M \varphi_{\varepsilon}(t)\sqrt{-1}\partial\bar{\partial}\frac{\partial\eta}{\partial t}\\\nonumber
&\xrightarrow{\varepsilon_{i}\rightarrow 0}&\int_M \omega_0\wedge\frac{\partial\eta}{\partial t}+\int_M \varphi(t)\sqrt{-1}\partial\bar{\partial}\frac{\partial\eta}{\partial t}\\\nonumber
&=&\int_M \omega_0\wedge\frac{\partial\eta}{\partial t}+\int_M \sqrt{-1}\partial\bar{\partial}\varphi(t)\frac{\partial\eta}{\partial t}\\
&=&\int_M \omega(t)\wedge\frac{\partial\eta}{\partial t}.
\end{eqnarray}

On the other hand, $\varphi_{\varepsilon}(t)$ and $\frac{\partial\varphi_\varepsilon(t)}{\partial t}$ are uniformly bounded on $[\delta,T]\times M$, $\varphi(t)$ and $\frac{\partial\varphi(t)}{\partial t}$ are uniformly bounded on $[\delta,T]\times(M\setminus D)$, therefore
\begin{eqnarray}\label{201503203}\nonumber\frac{\partial}{\partial t}\int_M \omega_{\varphi_{\varepsilon}(t)}\wedge\eta
&=&\int_M \varphi_{\varepsilon}(t)\sqrt{-1}\partial\bar{\partial}\frac{\partial\eta}{\partial t}\\\nonumber
&\ &+\int_M \frac{\partial\varphi_{\varepsilon}(t)}{\partial t} \sqrt{-1}\partial\bar{\partial}\eta+\int_M \omega_0\wedge\frac{\partial\eta}{\partial t}\end{eqnarray}
\begin{eqnarray}\nonumber
&\xrightarrow{\varepsilon\rightarrow 0}&
\int_M \varphi(t)\sqrt{-1}\partial\bar{\partial}\frac{\partial\eta}{\partial t}\\\nonumber
&\ &+\int_M \frac{\partial\varphi}{\partial t} \sqrt{-1}\partial\bar{\partial}\eta+\int_M \omega_0\wedge\frac{\partial\eta}{\partial t}\\
&=&\frac{\partial}{\partial t}\int_M \omega(t)\wedge\eta.
\end{eqnarray}
Combining equality
$$\frac{\partial}{\partial t}\int_M \omega_{\varphi_{\varepsilon}(t)}\wedge\eta=\int_{M}\frac{\partial\omega_{\varphi_{\varepsilon}(t)}}{\partial t}\wedge\eta+\int_M \omega_{\varphi_{\varepsilon}(t)}\wedge\frac{\partial\eta}{\partial t}$$
with equalities $(\ref{201503201})$-$(\ref{201503203})$, on $[\delta,T]$, we have
\begin{eqnarray}\label{201503204}\nonumber\frac{\partial}{\partial t}\int_M \omega(t)\wedge\eta&=&\int_{M} \big(-Ric (\omega(t))+\beta \omega(t) +2\pi (1-\beta ) [D]\big)\wedge\eta\\
&\ &+\int_M \omega(t)\wedge\frac{\partial\eta}{\partial t}.
\end{eqnarray}
Integrating form $0$ to $\infty$ on both sides,
\begin{eqnarray*}\int_{M\times(0,\infty)} \frac{\partial \omega(t)}{\partial t}\wedge\eta~dt
&=&-\int_{M\times(0,\infty)} \omega(t)\wedge\frac{\partial\eta}{\partial t}~dt=-\int_{0}^{\infty}\int_M \omega(t)\wedge\frac{\partial\eta}{\partial t}~dt\\
&=&\int_{0}^{\infty}\int_{M} \big(-Ric (\omega(t))+\beta \omega(t) +2\pi (1-\beta ) [D]\big)\wedge\eta~dt\\
&=&\int_{M\times(0,\infty)} \big(-Ric (\omega(t))+\beta \omega(t) +2\pi (1-\beta ) [D]\big)\wedge\eta~dt.
\end{eqnarray*}
By the arbitrariness of $\eta$, we prove that $\omega(t)$ satisfies the conical K\"ahler-Ricci flow $(\ref{TCKRF1})$ in the sense of currents on $(0,\infty)\times M$.  \QEDB

\medskip

Now we are ready to prove the uniqueness of the parabolic Monge-Amp\`ere equation $(\ref{223})$ starting with $\varphi_0\in\mathcal{E}_{p}(M,\omega_{0})$ for some $p>1$.

\begin{thm}\label{228} Let $\varphi_i(t)\in C^0\big([0,\infty)\times M\big)\bigcap C^\infty\big((0,\infty)\times(M\setminus D)\big)$ $(i=1,2)$ be two  long-time solutions to the parabolic Monge-Amp\`ere equation
\begin{eqnarray}
  \frac{\partial \varphi_i(t)}{\partial t}=\log\frac{(\omega_{0}+\sqrt{-1}\partial\bar{\partial}\varphi_i(t))^{n}}{\omega_{0}^{n}}+F_{0}+\beta\varphi_i(t)
  +\log|s|_{h}^{2(1-\beta)}
\end{eqnarray}
on $(0,\infty)\times (M\setminus D)$. If $\varphi_i$ $(i=1,2)$ satisfy
\begin{itemize}
  \item For any $0<\delta<T<\infty$, there exists uniform constant $C$ such that
\begin{eqnarray*}C^{-1}\omega_\beta\leq \omega_{0}+\sqrt{-1}\partial\bar{\partial}\varphi_i(t)\leq C\omega_\beta
\end{eqnarray*}
on $[\delta,T]\times (M\setminus D)$;
  \item On $[\delta, T]$, there exist constant $\alpha>0$ and $C^{*}$ such that $\varphi_i(t)$ is $C^{\alpha}$ on $M$ with respect to $\omega_{0}$ and $\| \frac{\partial\varphi_i(t)}{\partial t}\|_{L^{\infty}(M\setminus D)}\leqslant C^{*}$;
  \item $\lim\limits_{t\rightarrow0^{+}}\|\varphi_i(t)-\varphi_{0}\|_{L^{\infty}(M)}=0$.
  \end{itemize}
Then $\varphi_1=\varphi_2$.
\end{thm}

{\bf Proof:}\ \ We apply Jeffres' trick \cite{TJEF} in the parabolic case. For any $0<t_1<T<\infty$ and $a>0$. Let $\phi_1(t)=\varphi_1(t)+a|s|_h^{2q}$, where $0<q<1$ is determined later. The evolution of $\phi_1$ is
\begin{eqnarray*}
  \frac{\partial \phi_1(t)}{\partial t}=\log\frac{(\omega_{0}+\sqrt{-1}\partial\bar{\partial}\varphi_1(t))^{n}}{\omega_{0}^{n}}+F_{0}+\beta\phi_1(t)
  -a\beta|s|_h^{2q}+\log|s|_{h}^{2(1-\beta)}.
\end{eqnarray*}
Denote $\psi(t)=\phi_1(t)-\varphi_2(t)$ and $\hat{\Delta}=\int_0^1 g_{s\varphi_1+(1-s)\varphi_2}^{i\bar{j}}\frac{\partial^2}{\partial z^i\partial\bar{z}^j}ds$, $\psi(t)$ evolves along the following equation
\begin{eqnarray*}
  \frac{\partial \psi(t)}{\partial t}=\hat{\Delta}\psi(t)-a\hat{\Delta}|s|_h^{2q}+\beta\psi(t)
  -a\beta|s|_h^{2q}.
\end{eqnarray*}
By the equivalence of the metrics and the equation
\begin{eqnarray*}
\sqrt{-1}\partial\overline{\partial}|s|_h^{2q}=q^2|s|_h^{2q}\sqrt{-1}\partial\log|s|_h^{2}\wedge\overline{\partial}\log|s|_h^{2}
+q|s|_h^{2q}\sqrt{-1}\partial\overline{\partial}\log|s|_h^{2},
\end{eqnarray*}
we obtain the estimate
\begin{eqnarray*}
\hat{\Delta}|s|_h^{2q}&\geq &q|s|_h^{2q}g_{s\varphi_1+(1-s)\varphi_2}^{i\bar{j}}(\frac{\partial^2}{\partial z^i\partial\bar{z}^j}\log|s|_h^{2})\\
&=&-q|s|_h^{2q}g_{s\varphi_1+(1-s)\varphi_2}^{i\bar{j}}g_{0,i\bar{j}}\\
&\geq&-Cq|s|_h^{2q}g_{\beta}^{i\bar{j}}g_{0,i\bar{j}}\\
&\geq&-C
\end{eqnarray*}
on $M\setminus D$, where constant $C$ independent of $a$, and we apply the fact that $\omega_\beta\geq\gamma\omega_0$ on $M\setminus D$ for some constant $\gamma$. Then we obtain
\begin{eqnarray*}
  \frac{\partial \psi(t)}{\partial t}\leq\hat{\Delta}\psi(t)+\beta\psi(t)+aC.
\end{eqnarray*}

Let $\tilde{\psi}=e^{-\beta (t-t_1)}\psi+\frac{aC}{\beta}e^{-\beta (t-t_1)}-\epsilon (t-t_1)$. By choosing suitable $0<q<1$, we can  assume that the space maximum of $\tilde{\psi}$ on $[t_1,T]\times M$ is attained away from $D$. Let $(t_0,x_0)$ be the maximum point. If $t_0>t_1$, by the maximum principle, at $(t_0,x_0)$, we have
\begin{eqnarray*}
 0\leq (\frac{\partial }{\partial t}-\hat{\Delta})\tilde{\psi}(t)\leq -\epsilon,
\end{eqnarray*}
which is impossible, hence $t_0=t_1$. Then for $(t,x)\in [t_1,T]\times M$, we obtain
\begin{eqnarray*}
\psi(t,x)&\leq& e^{\beta T}\|\varphi_1(t_1,x)-\varphi_2(t_1,x)\|_{L^\infty(M)}\\
&\ &+aCe^{\beta T}+\epsilon Te^{\beta T}
\end{eqnarray*}
Let $a\rightarrow0$ and then $t_1\rightarrow0^+$, we get
\begin{eqnarray*}
\varphi_1(t)-\varphi_2(t)\leq \epsilon Te^{\beta T}.
\end{eqnarray*}
It shows that $\varphi_1(t)\leq\varphi_2(t)$ after we let $\epsilon\rightarrow0$. By the same reason we have $\varphi_2(t)\leq\varphi_1(t)$, then we prove that $\varphi_1(t)=\varphi_2(t)$.\QEDB

\begin{lem} \label{20160117}Assume that on $(0,\infty)\times(M\setminus D)$, $f(t,x)$ is a smooth function and satisfis $\sqrt{-1}\partial\bar{\partial}f(t,x)=0$. If
\begin{eqnarray}
\|f(t,x)\|_{L^{\infty}([\delta,T]\times (M\setminus D))}\leqslant C_{\delta, T},
\end{eqnarray}
where $0<\delta<T<\infty$. Then $f(t,x)=f(t)$, i.e. function $f(t,x)$ depends only on $t$ on $(0,\infty)\times(M\setminus D)$.
\end{lem}

{\bf Proof:}\ \ We prove this lemma by the logarithmic cutoff trick. Fix a cutoff function $\eta:[0,\infty)\rightarrow [0,1]$ such that $\eta(s)=1$ for $s<1$ and $\eta(s)=0$ for $s>2$. We define
\begin{eqnarray}\gamma(x)=\eta(\frac{\log r}{\log \varepsilon}),
\end{eqnarray}
where $r$ is the the distance function from the divisor. Then $\gamma=0$ for $r<\varepsilon^{2}$ and $\gamma=1$ for $r>\varepsilon$. Straightforward calculations shows that
\begin{eqnarray*}\int_{M} |\partial \gamma|_{g_{0}}^{2} dV_{0}&=&\int_{M}(\eta^{'})^{2} \frac{1}{r^{2}(\log\varepsilon)^{2}}|\partial r|_{g_{0}}^{2}dV_{0}\\
&\leqslant&C\int_{\varepsilon^{2}}^{\varepsilon}\frac{1}{r^{2}(\log\varepsilon)^{2}} rdr=-C\frac{1}{\log\varepsilon},
\end{eqnarray*}
so we have $\int_{M} |\partial \gamma|_{g_{0}}^{2} dV_{0}\rightarrow0$ as $\varepsilon\rightarrow0$. On the other hand, when $t\in[\delta, T]$,
\begin{eqnarray*}
0&=&\int_{M} div_{g_{0}}\big(\gamma^{2}f(\overline{\nabla}_{g_{0}}f)\big)dV_{0}\\
&=&2\int_{M} \gamma f<\partial \gamma,\overline{\partial }f>dV_{0}+\int_{M} \gamma^{2}|\partial f|_{g_{0}}^{2}dV_{0}\\
&\geqslant&\frac{1}{2}\int_{M} \gamma^{2}|\partial f|_{g_{0}}^{2}dV_{0}-2\int_{M} f^{2}|\partial \gamma|_{g_{0}}^{2}dV_{0},
\end{eqnarray*}
which implies
\begin{eqnarray}
\int_{M} \gamma^{2}|\partial f|_{g_{0}}^{2}dV_{0}\leqslant C\int_{M} |\partial \gamma|_{g_{0}}^{2}dV_{0},
\end{eqnarray}
where constant $C$ depends only on $\|f\|_{L^{\infty}([\delta,T]\times (M\setminus D))}$. Passing to the limit and then combining with the arbitrary choice of $\delta$ and $T$, we obtain
\begin{eqnarray}
\int_{M} |\partial f|_{g_{0}}^{2}dV_{0}=0\ \ \ \ \ on\ \ \ (0,\infty).
\end{eqnarray}
Hence function $f(t,x)$ depends only on $t$. \QEDB

\begin{thm}\label{229} $\omega_{\varphi(t)}=\omega_{0}+\sqrt{-1}\partial\bar{\partial}\varphi(t)$ is the unique long-time solution to the conical K\"ahler-Ricci flow $(\ref{TCKRF1})$.
\end{thm}

{\bf Proof:}\ \ Suppose there is another solution $\omega_{\phi(t)}=\omega_0+\sqrt{-1}\partial\bar{\partial}\phi(t)$ to the conical K\"ahlre-Ricci flow $(\ref{TCKRF1})$. It follows from Lemma $\ref{20160117}$ that
\begin{eqnarray}\label{20150618}
\frac{\partial\phi(t)}{\partial t}=\log\frac{(\omega_{0}+\sqrt{-1}\partial\bar{\partial}\phi(t))^{n}}{\omega_{0}^{n}}+F_{0}+\beta\phi(t)
  +\log|s|_{h}^{2(1-\beta)}+f(t)
\end{eqnarray}
on $(0,\infty)\times (M\setminus D)$ for a smooth function $f(t)$ defined on $(0,\infty)$, and $\phi(t)\in C^0\big([0,\infty)\times M\big)\bigcap C^\infty\big((0,\infty)\times(M\setminus D)\big)$ satisfies
\begin{itemize}
  \item For any $0<\delta<T<\infty$, there exists uniform constant $C$ such that
\begin{eqnarray*}C^{-1}\omega_\beta\leq \omega_{0}+\sqrt{-1}\partial\bar{\partial}\phi(t)\leq C\omega_\beta\ \ \ \ on\ \ \ \ [\delta,T]\times (M\setminus D);
\end{eqnarray*}

  \item On $[\delta, T]$, there exist constant $\alpha>0$ and $C$ such that $\phi(t)$ is $C^{\alpha}$ on $M$ with respect to $\omega_{0}$ and $\| \frac{\partial\phi(t)}{\partial t}\|_{L^{\infty}(M\setminus D)}\leqslant C$;
  \item $\lim\limits_{t\rightarrow0^{+}}\|\phi(t)-\varphi_{0}\|_{L^{\infty}(M)}=0$.
  \end{itemize}

For any $0<t_1<T<\infty$ and $a>0$. Let $\psi(t)=\phi(t)+a|s|_h^{2q}-\varphi(t)$, where $0<q<1$ is determined later. Then
\begin{eqnarray*}
  \frac{\partial \psi(t)}{\partial t}=\hat{\Delta}\psi(t)-a\hat{\Delta}|s|_h^{2q}+\beta\psi(t)
  -a\beta|s|_h^{2q}+f(t).
\end{eqnarray*}
By the same arguments as that in the proof of Proposition \ref{228}, for any $(t,x)\in [t_1,T]\times M$, we have
\begin{eqnarray*}
\psi(t,x)&\leq& e^{\beta (t-t_1)}\|\phi(t_1,x)-\varphi(t_1,x)\|_{L^\infty(M)}\\
&\ &+aCe^{\beta (t-t_1)}+\epsilon (t-t_1)e^{\beta (t-t_1)}\\
&\ &+e^{\beta (t-t_1)}\int_{t_1}^te^{-\beta (s-t_1)}f(s)ds
\end{eqnarray*}
Let $a\rightarrow0$, we obtain
\begin{eqnarray*}
\phi(t)-\varphi(t)&\leq& e^{\beta (t-t_1)}\|\phi(t_1,x)-\varphi(t_1,x)\|_{L^\infty(M)}\\
&\ &+\epsilon (t-t_1)e^{\beta (t-t_1)}+e^{\beta (t-t_1)}\int_{t_1}^te^{-\beta (s-t_1)}f(s)ds.
\end{eqnarray*}
By the similar arguments, we can obtain
\begin{eqnarray*}
\varphi(t)-\phi(t)&\leq& e^{\beta (t-t_1)}\|\phi(t_1,x)-\varphi(t_1,x)\|_{L^\infty(M)}\\
&\ &+\epsilon (t-t_1)e^{\beta (t-t_1)}-e^{\beta (t-t_1)}\int_{t_1}^te^{-\beta (s-t_1)}f(s)ds.
\end{eqnarray*}
Therefore, for any $t>t_1>0$, we have
\begin{eqnarray*}
\inf\limits_{M}(\phi(t)-\varphi(t))&\geq& \sup\limits_{M}(\phi(t)-\varphi(t))-2e^{\beta (t-t_1)}\|\phi(t_1,x)-\varphi(t_1,x)\|_{L^\infty(M)}\\
&\ &-2\epsilon (t-t_1)e^{\beta (t-t_1)}\\
&\geq& \sup\limits_{M}(\phi(t)-\varphi(t))-2e^{\beta T}\|\phi(t_1,x)-\varphi(t_1,x)\|_{L^\infty(M)}\\
&\ &-2\epsilon Te^{\beta T}.
\end{eqnarray*}
Let $t_1\rightarrow0^+$ and then $\epsilon\rightarrow0$, we conclude that $\phi(t)=\varphi(t)+e^{\beta t}\int_{0}^te^{-\beta s}f(s)ds$. Then
$\omega_{\phi(t)}=\omega_{\varphi(t)}$ on $(0,\infty)\times (M\setminus D)$. \QEDB

\section{The convergence of the conical K\"ahler-Ricci flow with weak initial data}
\setcounter{equation}{0}

In this section, we study the convergence of the conical K\"ahler-Ricci flow (\ref{TCKRF1}) on Fano manifold with positive twisted first Chern class. Our discussion is very similar as that in \cite{JWLXZ}, but we need new arguments on estimates of the twisted Ricci potential $u_{\varepsilon}(t)$ and the term $|\dot{\varphi}_{\varepsilon}|$ when we handle the weak initial data case.

  Without loss of generality, we assume $\lambda=1$ (i.e. $\mu=\beta$). We first prove the uniform Perelman's estimates along the twisted K\"ahler-Ricci flow
\begin{eqnarray}\label{TKRF3}
\begin{cases}
  \frac{\partial \omega_{\varepsilon}(t)}{\partial t}=-Ric(\omega_{\varepsilon}(t))+\beta\omega_{\varepsilon}(t)+\theta_{\varepsilon}.\\
  \\
  \omega_{\varepsilon}(t)|_{t=0}=\omega_{\varphi_{0}}\\
  \end{cases}
  \end{eqnarray}
By the same argument as Proposition $4.1$ in \cite{JWLXZ}, we have

\medskip

\begin{pro}\label{1.8.1} $t^2(R(g_{\varepsilon,j}(t))-tr_{g_{\varepsilon,j}(t)}\theta_{\varepsilon})$ is uniformly bounded from below along the twisted K\"ahler-Ricci flow $(\ref{TKRF2})$, i.e. there exists a uniform constant $C$, such that
\begin{equation}\label{3.22.9}t^2(R(g_{\varepsilon,j}(t))-tr_{g_{\varepsilon,j}(t)}\theta_{\varepsilon})\geq-C \end{equation}
for any $t\geqslant0$, $j\in\mathbb{N}^{+}$ and $\varepsilon>0$, while the constant $C$ only depends on $\beta$ and $n$. In particular,
\begin{equation}\label{3.22.911}R(g_{\varepsilon,j}(t))-tr_{g_{\varepsilon,j}(t)}\theta_{\varepsilon}\geq-C \end{equation}
when $t\geq\frac{1}{2}$.
\end{pro}

\begin{rem}\label{232}By Proposition \ref{210}, we know that there exists constant $C$ only depending on $\beta$ and $n$, such that
\begin{equation}\label{ccc1}
R(g_{\varepsilon}(t))-tr_{g_{\varepsilon}(t)}\theta_{\varepsilon}\geq-C \end{equation}
along the twisted K\"ahler-Ricci flow $(\ref{TKRF3})$ for any $\varepsilon>0$ when $t\geq\frac{1}{2}$.
\end{rem}

Straightforward calculation shows that the twisted Ricci potential $u_\varepsilon(t)$ with respect to $\omega_\varepsilon(t)$ at $t=\frac{1}{2}$ can be written as
\begin{eqnarray}
u_\varepsilon(\frac{1}{2})=\log\frac{\omega_\varepsilon^n(\frac{1}{2})(\varepsilon^2+|s|_h^2)^{1-\beta}}{\omega_0^n}+F_0
+\beta\varphi_\varepsilon(\frac{1}{2})+C_{\varepsilon,\frac{1}{2}},
\end{eqnarray}
where $C_{\varepsilon,\frac{1}{2}}$ is a normalization constant such that $\frac{1}{V}\int_M e^{-u_\varepsilon(\frac{1}{2})}dV_{\varepsilon,\frac{1}{2}}=1$. By $(\ref{219})$ and
$(\ref{220})$, we conclude that $C_{\varepsilon,\frac{1}{2}}$ and $u_{\varepsilon}(\frac{1}{2})$ are uniformly bounded. Let $a_{\varepsilon}(t)=\frac{\beta}{V}\int_{M}u_{\varepsilon}(t)e^{-u_{\varepsilon}(t)}dV_{\varepsilon,t}$, then by Lemma $4.4$ in \cite{JWLXZ}, we have

\begin{lem}\label{1.8.4} There exists a uniform constant $C$, such that
\begin{equation}\label{3.22.8}|a_{\varepsilon}(t)|\leq C\end{equation}
for any $t\geq\frac{1}{2}$ and $\varepsilon>0$.
\end{lem}

Now we consider the twisted K\"ahler-Ricci flows $(\ref{TKRF3})$ starting at $t=\frac{1}{2}$. Using the estimates (\ref{ccc1}), (\ref{3.22.8}) and following the arguments in section $4$ of \cite{JWLXZ}, we have the following uniform Perelman's estimates.

\begin{thm}\label{1.8.2} Let $g_{\varepsilon}(t)$ be a solution of the twisted K$\ddot{a}$hler Ricci flow, i.e. the corresponding form $\omega_{\varepsilon}(t)$ satisfies the equation $(\ref{TKRF3})$ with initial metric $\omega_{\varphi_{0}}$, $u_{\varepsilon}(t)\in C^{\infty}((0,\infty)\times M)$ is the twisted Ricci potential satisfying
\begin{equation}\label{4.1}-Ric(\omega_{\varepsilon}(t))+\beta\omega_{\varepsilon}(t)+\theta_{\varepsilon}=\sqrt{-1}\partial\bar{\partial}u_{\varepsilon}(t)\end{equation}
and $\frac{1}{V}\int_{M}e^{-u_{\varepsilon}(t)} dV_{\varepsilon, t}=1$, where $\theta_{\varepsilon}=(1-\beta)(\omega_{0}+\sqrt{-1}\partial\overline{\partial}\log(\varepsilon^{2}+|s|_{h}^{2}))$. Then for any $\beta\in(0,1)$, there exists a uniform constant $C$, such that
\begin{eqnarray}\label{p0}
|R(g_{\varepsilon}(t))-tr_{g_{\varepsilon}(t)}\theta_{\varepsilon}|&\leq& C,\\
\label{p1}\|u_{\varepsilon}(t)\|_{C^{1}(g_{\varepsilon}(t))}&\leq& C,\\
\label{p2}diam(M,g_{\varepsilon}(t))&\leq&C\end{eqnarray}
hold for any $t\geq1$ and $\varepsilon>0$, where $R(g_{\varepsilon}(t))-tr_{g_{\varepsilon}(t)}\theta_{\varepsilon}$ and $diam(M,g_{\varepsilon}(t))$ are the twisted scalar curvature and diameter of the manifold respectively with respect to the metric $g_{\varepsilon}(t)$.
\end{thm}

\medskip

If $\varphi_\varepsilon(t)$ is a solution to the Monge-Amp\`ere equation
\begin{eqnarray}\label{330}
\frac{\partial \varphi_{\varepsilon}(t)}{\partial t}=\log\frac{(\omega_{0}+\sqrt{-1}\partial\bar{\partial}\varphi_{\varepsilon}(t))^{n}}{\omega_{0}^{n}}+F_{0}+\beta\varphi_{\varepsilon}(t)
  +\log(\varepsilon^{2}+|s|_{h}^{2})^{1-\beta}
\end{eqnarray}
on $(0,\infty)\times M$ with initial value $\varphi_\varepsilon(0)=\varphi_0$, it is obvious that $\phi_\varepsilon(t)=\varphi_\varepsilon(t)+Ce^{\beta t}$ is a solution to equation $(\ref{330})$ with initial value $\phi_\varepsilon(0)=\varphi_0+C$. At the same time, $\omega_{\phi_\varepsilon(t)}=\omega_0+\sqrt{-1}\partial\bar{\partial}\phi_\varepsilon(t)$ is also a solution to the twisted K\"ahler-Ricci flow $(\ref{TKRF3})$ with initial value $\omega_{\varphi_0}$.

From $(\ref{219})$, we know that $\varphi_{\varepsilon}(t)$ is uniformly bounded on $[0,T]\times M$ by a constant $C$ which depends only on $\|\varphi_0\|_{L^\infty(M)}$, $\beta$ and $T$. Now, we consider the solution $\psi_\varepsilon(t)=\varphi_\varepsilon(t)+\tilde{C}_{\varepsilon,1}e^{\beta t}$ to the equation
\begin{eqnarray}\label{331}
\begin{cases}
  \frac{\partial \psi_{\varepsilon}(t)}{\partial t}=\log\frac{(\omega_{0}+\sqrt{-1}\partial\bar{\partial}\psi_{\varepsilon}(t))^{n}}{\omega_{0}^{n}}+F_{0}+\beta\psi_{\varepsilon}(t)
  +\log(\varepsilon^{2}+|s|_{h}^{2})^{1-\beta}\\
  \ \ \ \ \ \ \ \ \ \ \ \ \ \ \ \ \ \ \ \ \ \ \ \ \ \ \ \ \ \ \ \ \ \ \ \ \ \ \ \ \ \ \ \ \ \ \ \ \ \ \ \ \ \ \ \ \ \ \ \ \ \ on\ \ (0,\infty)\times M,\\
  \psi_{\varepsilon}(0)=\varphi_{0}+\tilde{C}_{\varepsilon,1}\\
  \end{cases}
\end{eqnarray}
where
\begin{eqnarray*}
\tilde{C}_{\varepsilon,1}=e^{-\beta}\frac{1}{\beta}\Big(\int_{1}^{+\infty}e^{-\beta t}\|\nabla u_{\varepsilon}(t)\|^{2}_{L^{2}}dt-\frac{1}{V}\int_{M}\big(F_{\varepsilon,1}+\beta\varphi_\varepsilon(1)\big)dV_{\varepsilon,1}\Big),
\end{eqnarray*} $F_{\varepsilon,1}=F_{0}+\log(\frac{\omega_{\varepsilon}(1)^{n}}{\omega_{0}^{n}}\cdot(\varepsilon^{2}+|s|_{h}^{2})^{1-\beta})$ and $dV_{\varepsilon,1}=\frac{\omega_\varepsilon^n(1)}{n!}$. By $(\ref{219})$, $(\ref{220})$ and the above uniform Perelman's estimates $(\ref{p1})$,  we know that the constant $\tilde{C}_{\varepsilon,1}$ is well-defined and uniformly bounded. Straightforward calculation shows that the twisted Ricci potential $u_{\varepsilon}(1)$ with respect to $\omega_\varepsilon(1)$ can be written as
\begin{eqnarray}\label{332}
u_\varepsilon(1)=\log\frac{\omega_\varepsilon^n(1)(\varepsilon^2+|s|_h^2)^{1-\beta}}{\omega_0^n}+F_0
+\beta\varphi_\varepsilon(1)+C_{\varepsilon,1},
\end{eqnarray}
where $C_{\varepsilon,1}$ is a normalization constant such that $\frac{1}{V}\int_M e^{-u_\varepsilon(1)}dV_{\varepsilon,1}=1$. Then
\begin{eqnarray}
C_{\varepsilon,1}=\log\Big(\frac{1}{V}\int_M e^{-F_0-\beta\varphi_\varepsilon(1)}\frac{dV_{0}}{(\varepsilon^2+|s|_h^2)^{1-\beta}}\Big).
\end{eqnarray}
By $(\ref{219})$ and
$(\ref{220})$, we conclude that $C_{\varepsilon,1}$ and $u_{\varepsilon}(1)$ are uniformly bounded.

Let $u_\varepsilon(t)= \dot{\psi}_\varepsilon(t)+c_\varepsilon(t)$. By equation $(\ref{331})$ and equality $(\ref{332})$, we have
\begin{eqnarray}\label{333}
c_{\varepsilon}(1)=C_{\varepsilon,1}-\beta e^\beta\tilde{C}_{\varepsilon,1}.
\end{eqnarray}

\begin{pro}\label{1.8.5.3}There exists a uniform constant $C$ such that
$$\|\dot{\psi}_{\varepsilon}(t)\|_{C^{0}}\leq C$$
for any $\varepsilon>0$ and $t\geq1$.
\end{pro}

{\bf  Proof:}\ \ As in \cite{PSS}, when $t\geqslant1$, we let
\begin{eqnarray}\label{3.22.32}\alpha_{\varepsilon}(t)=\frac{1}{V}\int_{M}\dot{\psi}_{\varepsilon}(t)dV_{\varepsilon,t}
=\frac{1}{V}\int_{M}u_{\varepsilon}(t)dV_{\varepsilon,t}-c_{\varepsilon}(t).
\end{eqnarray}
Through computing, we have
\begin{eqnarray}\frac{d}{dt}\alpha_{\varepsilon}(t)&=&\beta\alpha_{\varepsilon}(t)-\|\nabla\dot{\psi}_{\varepsilon}\|^{2}_{L^{2}},
\end{eqnarray}
\begin{eqnarray}
\label{3.22.33}\nonumber e^{-\beta (t-1)}\alpha_{\varepsilon}(t)&=&\alpha_{\varepsilon}(1)-\int_{1}^{t}e^{-\beta (s-1)}\|\nabla\dot{\psi}_{\varepsilon}\|^{2}_{L^{2}}ds\\
&=&\frac{1}{V}\int_{M}u_{\varepsilon}(1)dV_{\varepsilon,1}-c_{\varepsilon}(1)-\int_{1}^{t}e^{-\beta (s-1)}\|\nabla\dot{\psi}_{\varepsilon}\|^{2}_{L^{2}}ds.
\end{eqnarray}
Putting $(\ref{332})$ and $(\ref{333})$ into $(\ref{3.22.33})$, we have
\begin{eqnarray}e^{-\beta (t-1)}\alpha_{\varepsilon}(t)=\int_{t}^{+\infty}e^{-\beta (s-1)}\|\nabla\dot{\psi}_{\varepsilon}\|^{2}_{L^{2}}ds.
\end{eqnarray}
By Theorem \ref{1.8.2}, we conclude that
\begin{eqnarray}\label{3.22.34}0\leq\alpha_{\varepsilon}(t)&=&\int_{t}^{+\infty}e^{\beta(t-s)}\|\nabla\dot{\psi}_{\varepsilon}\|^{2}_{L^{2}}ds\leq C.
\end{eqnarray}
Then we conclude that $\dot{\psi}_{\varepsilon}(t)$ is uniformly bounded by the uniform Perelman's estimates when $t\geqslant1$.\QEDB

\medskip

We recall Aubin's functionals, Ding's functional and the twisted Mabuchi $\mathcal{K}$-energy functional.
\begin{eqnarray}\label{3.22.35}
I_{\omega_0}(\phi)&=&\frac{1}{V}\int_M\phi(dV_{0}-dV_{\phi}), \\ \nonumber
J_{\omega_0}(\phi)&=&\frac{1}{V}\int_0^1\int_M\dot{\phi}_t(dV_{0}-dV_{\phi_{t}})dt,
\end{eqnarray}
where $\phi_t$ is a path with $\phi_0=c$, $\phi_1=\phi$.
\begin{eqnarray}\label{3.22.360}
F^0_{\omega_0}(\phi)&=&J_{\omega_0}(\phi)-\frac{1}{V}\int_M\phi dV_0,\\
\label{3.22.361}F_{\omega_0, \theta}(\phi)&=&J_{\omega_0}(\phi)-\frac{1}{V}\int_M\phi dV_0-\frac{1}{\beta}\log(\frac{1}{V}\int_Me^{-u_{\omega_{0}}-\beta\phi}dV_0),\\
\label{3.22.36'}\mathcal{M}_{ \omega_{0},\ \theta}(\phi)&=&-\beta(I_{\omega_{0}}(\phi)-J_{\omega_{0}}(\phi))-\frac{1}{V}\int_{M}u_{\omega_{0}}(dV_{0}-dV_{\phi})\\ \nonumber
&\ &+\frac{1}{V}\int_{M}\log\frac{\omega_{\phi}^{n}}{\omega_{0}^{n}}dV_{\phi},
\end{eqnarray}
where $u_{\omega_{0}}$ is the twisted Ricci potential of $\omega_0$, $i.e.$ $-Ric(\omega_{0})+\beta\omega_{0}+\theta=\sqrt{-1}\partial\bar{\partial} u_{\omega_{0}}$ and $\frac{1}{V}\int_{M}e^{-u_{\omega_{0}}}dV_{\omega_{0}}=1$.

\begin{pro}\label{1.8.5.4} For any $t\geq1$, the solution $\psi_{\varepsilon}(t)$ to equation $(\ref{331})$ satisfies:
\begin{eqnarray*}&(i)& \ \ \ \mathcal{M}_{ \omega_{0},\ \theta_{\varepsilon} }(\psi_{\varepsilon}(t))-\beta F^0_{\omega_0}(\psi_{\varepsilon}(t))-\frac{1}{V}\int_M\dot{\psi}_{\varepsilon}(t) dV_{\varepsilon, t}=C_{\varepsilon},\\
&(ii)& \ \ \mathcal{M}_{ \omega_{0},\ \theta_{\varepsilon}}(\psi_{\varepsilon}(1))\ is\ uniformly\ bounded,\end{eqnarray*}
where $C_{\varepsilon}$ in $(i)$ can be bounded by a uniform constant $C$.
\end{pro}

{\bf  Proof:}\ \ Following the argument in \cite{JWLXZ}, since
\begin{eqnarray}\label{3.22.37}\frac{d}{dt}(\mathcal{M}_{\omega_{0},\ \theta_{\varepsilon} }(\psi_{\varepsilon}(t))-\beta F^0_{\omega_0}(\psi_{\varepsilon}(t))-\frac{1}{V}\int_M\dot{\psi}_{\varepsilon}(t) dV_{\varepsilon, t})=0,
\end{eqnarray}
we obtain that
\begin{eqnarray*}&\ &\mathcal{M}_{\omega_{0},\ \theta_{\varepsilon} }(\psi_{\varepsilon}(t))-\beta F^0_{\omega_0}(\psi_{\varepsilon}(t))-\frac{1}{V}\int_M\dot{\psi}_{\varepsilon}(t) dV_{\varepsilon, t}\\
&=&\mathcal{M}_{\omega_{1},\ \theta_{\varepsilon} }(\psi_{\varepsilon}(1))-\beta F^0_{\omega_0}(\psi_{\varepsilon}(1))-\frac{1}{V}\int_M\dot{\psi}_{\varepsilon}(1) dV_{\varepsilon,1}\\
&=&\frac{1}{V}\int_{M}\log\frac{\omega_{\varepsilon}^{n}(1)(|s|^{2}_{h}+\varepsilon^{2})^{1-\beta}}{e^{-F_{0}}\omega_{0}^{n}}dV_{\varepsilon,1}
+\frac{\beta }{V}\int_{M}\psi_{\varepsilon}(1)dV_{\varepsilon,1}\\
&\ &-\frac{1}{V}\int_{M}F_{0}+\log(|s|^{2}_{h}+\varepsilon^{2})^{1-\beta}dV_{0}-\frac{1}{V}\int_M\dot{\psi}_{\varepsilon}(1) dV_{\varepsilon,1}.
\end{eqnarray*}
where the last equality can be bounded by a uniform constant. This gives a proof of $(i)$.
Furthermore, by the definition of $\mathcal{M}_{\omega_{0},\ \theta_{\varepsilon} }$, we have
\begin{eqnarray*}\mathcal{M}_{\omega_{0},\ \theta_{\varepsilon} }(\psi_{\varepsilon}(1))
&=&\frac{1}{V}\int_{M}\log\frac{\omega_{\varepsilon}^{n}(1)(|s|^{2}_{h}+\varepsilon^{2})^{1-\beta}}{e^{-F_{0}}\omega_{0}^{n}}dV_{\varepsilon,1}-\beta I_{\omega_{0}}(\psi_{\varepsilon}(1))+\beta J_{\omega_{0}}(\psi_{\varepsilon}(1))\\
&\ &-\frac{1}{V}\int_{M}F_{0}+\log(|s|^{2}_{h}+\varepsilon^{2})^{1-\beta}dV_{0}.
\end{eqnarray*}
Since $I_{\omega_{0}}(\psi_{\varepsilon}(1))$ is uniformly bounded and $\frac{1}{n}J_{\omega_0}\leq\frac{1}{n+1}I_{\omega_0}\leq J_{\omega_0}$, we prove $(ii)$.\QEDB

\medskip

Using Proposition \ref{1.8.5.3} and \ref{1.8.5.4}, by following the arguments in section $5$ of \cite{JWLXZ}, we obtain the following uniform $C^0$ estimate of $\psi_\varepsilon(t)$ along the equation $(\ref{331})$ under the assumption that the twisted Mabuchi $\mathcal{K}$-energy functional $\mathcal{M}_{ \omega_{0},\ \theta_{\varepsilon} }$ is uniformly proper on the space
\begin{eqnarray}\label{3.22.41}\mathcal{H}(\omega_0)=\{\phi\in C^{\infty}(M)|\ \omega_0+\sqrt{-1}\partial\bar{\partial}\phi>0\}.
\end{eqnarray}

\begin{thm}\label{1.8.5.6} Let $\psi_{\varepsilon}(t)$ be a solution of the flow $(\ref{331})$. If the twisted Mabuchi $\mathcal{K}$-energy functional $\mathcal{M}_{ \omega_{0},\ \theta_{\varepsilon} }$ is uniformly proper on $\mathcal{H}(\omega_0)$, i.e. there exists a uniform function $f$ such that
\begin{eqnarray}\label{3.22.42}\mathcal{M}_{ \omega_{0},\ \theta_{\varepsilon}}(\phi)\geq f(J_{\omega_0}(\phi))
\end{eqnarray}
for any $\varepsilon$ and $\phi\in\mathcal{H}(\omega_0)$, where $f(t):\mathbb{R}^+\rightarrow\mathbb{R}$ is some monotone increasing function satisfying $\lim\limits_{t\rightarrow+\infty}f(t)=+\infty$, then there exists a uniform constant $C$ such that for any $\varepsilon>0$ and $t\geq0$
\begin{eqnarray}\label{3.22.43}\|\psi_{\varepsilon}(t)\|_{C^{0}}\leq C.
\end{eqnarray}
\end{thm}

We study the flow $(\ref{TKRF3})$ start at $t=\frac{1}{2}$. We  obtain
$C^{-1}\omega_{\varepsilon}\leqslant \omega_{\varepsilon}(\frac{1}{2})\leqslant C\omega_{\varepsilon}$ on $M$ and
 $\|\psi_{\varepsilon}(\frac{1}{2})\|_{C^{k}(K)} \leqslant C_{k,K}$ on $K\subset\subset M\setminus D$
for some uniform constants $C$ and $C_{k,K}$ in section $3$, after getting the uniform bound of $\dot{\psi}_{\varepsilon}(t)$ and $\psi_{\varepsilon}(t)$, we can prove the uniform Laplacian $C^{2}$ estimates and local uniform estimates for any $t\geqslant1$ and $\varepsilon>0$ by the arguments in \cite{JWLXZ} (see Proposition $2.1$ and $2.3$ in \cite{JWLXZ}). In fact, we prove the following theorem.

\begin{thm}\label{1.8.5.6.1} Under the assumption in Theorem \ref{1.8.5.6}, for any $k\in \mathbb{N}^{+} $ and $K\subset\subset M\setminus D$, there exists constant $C_{k,K}$ depending only on $\|\varphi_0\|_{L^\infty(M)}$, $n$, $\beta$, $k$, $\omega_0$ and $dist_{\omega_{0}}(K,D)$, such that for any $\varepsilon>0$ and $t\geqslant1$, we have
\begin{eqnarray}\|\psi_{\varepsilon}(t)\|_{C^{k}(K)}\leq C_{k,K}.
\end{eqnarray}
\end{thm}

\medskip

Now we assume that there exists a conical K\"ahler-Einstein metric with cone angle $2\pi \beta $ along $D$. When $\lambda>0$ and there is no nontrivial holomorphic field which is tangent to $D$ along $D$, Tian-Zhu  \cite{GTXHZ05} obtained the following Moser-Trudinger type inequality
\begin{eqnarray}\label{20150611}F_{\omega_0,(1-\beta)D}(\phi)\geqslant\delta J_{\omega_{0}}(\phi)-C,\ \ \ \ \ \ \forall \phi\in\mathcal{H}(\omega_0)
\end{eqnarray}
for some constants $\delta$ and $C$, where
\begin{eqnarray*}
F_{\omega_{0},(1-\beta)D}(\phi)=J_{\omega_{0}}(\phi)-\frac{1}{V}\int_{M}\phi dV_{0}-\frac{1}{\beta}\log\big(\frac{1}{V}\int_M\frac{1}{|s|_{h}^{2(1-\beta)}}e^{-F_{0}-\beta\phi}dV_0\big)\end{eqnarray*}
is defined in \cite{GTXHZ05} (see also \cite{LS}).

\begin{rem}
When $\lambda\geqslant 1$, R. Berman \cite{RB}, Li-Sun \cite{LS} proved that there is no nontrivial holomorphic vector field on $M$ tangent to divisor $D$, and Li-Sun also proved that the existence of conical K\"ahler-Einstein metric can deduce the properness of the Log Mabuchi $\mathcal{K}$-energy functional (see also Song-Wang's results in \cite{SW}).
\end{rem}

By the definition of $F_{\omega_0,\theta_{\varepsilon}}$ and $F_{\omega_{0},(1-\beta)D}$, we have
\begin{eqnarray}\label{201506112}\nonumber F_{\omega_0,\theta_{\varepsilon}}(\phi)-F_{\omega_0,(1-\beta)D}(\phi)&=&\frac{1}{\beta}\log\Big(\frac{1}{V}\int_Me^{-F_0-C_{0}-\beta\phi}\frac{dV_{0}}{|s|_h^{2(1-\beta)}}\Big)\\
&\ &-\frac{1}{\beta}\log\Big(\frac{1}{V}\int_Me^{-F_0-C_{\varepsilon}-\beta\phi}\frac{dV_0}{(\varepsilon^2+|s|_h^2)^{(1-\beta)}}\Big)\\ \nonumber
&\geq& -C,
\end{eqnarray}
where $C_{0}$ and $C_{\varepsilon}$ are two normalized constants, and $C$ is a constant independent of $\varepsilon$. So the Ding's functional $F_{\omega_0, \theta_{\varepsilon}}$ is uniform proper. By the normalization and Jensen's inequality, we have
\begin{eqnarray}\label{201506111}\frac{1}{V}\int_M -u_{\omega_\phi} dV_{\phi}\leq\log\big(\frac{1}{V}\int_M e^{-u_{\omega_\phi}} dV_{\phi}\big)=0.
\end{eqnarray}
Then we have the following inequalities by $(\ref{3.22.361})$, $(\ref{3.22.36'})$ and $(\ref{201506111})$.
\begin{eqnarray}\label{201506113}\nonumber
\mathcal{M}_{ \omega_{0},\ \theta_{\varepsilon}}(\phi)&=&\beta F_{\omega_0,\theta_{\varepsilon}}(\phi)+\frac{1}{V}\int_M u_{\omega_\phi} dV_{\phi}-\frac{1}{V}\int_M u_{\omega_0} dV_{0}\\\nonumber
&\geq&\beta F_{\omega_0,\theta_{\varepsilon}}(\phi)-\frac{1}{V}\int_M F_0+C_{\varepsilon}+(1-\beta)\log(\varepsilon^2+|s|_h^2) dV_{0}\\
&\geq&\beta F_{\omega_0,\theta_{\varepsilon}}(\phi)-C,
\end{eqnarray}
where constant $C$ independent of $\varepsilon$. Hence we deduce the uniform properness of the twisted Mabuchi $\mathcal{K}$-energy functional by $(\ref{20150611})$, $(\ref{201506112})$ and $(\ref{201506113})$, i.e.
\begin{eqnarray}\label{201506114}{M}_{ \omega_{0},\ \theta_{\varepsilon}}(\phi)\geqslant C_{1} J_{\omega_{0}}(\phi)-C_{2},\ \ \ \ \ \ \forall \phi\in\mathcal{H}(\omega_0)
\end{eqnarray}
for some uniform constants $C_{1}$ and $C_{2}$. At the same time, we have the uniqueness theorem of conical K\"ahler-Einstein metric (proved by B. Berndtsson in \cite{BBERN})  under the assumption that there is no nontrivial holomorphic field which is tangent to $D$. Using the above $C^{0}$ estimate and the uniqueness theorem, we can apply the arguments in section $6$ of \cite{JWLXZ} to obtain the convergence  result  of the conical K\"ahler-Ricci flow $(\ref{TCKRF1})$, i.e. Theorem 1.3.

\hspace{1.4cm}

\end{document}